\documentclass[reqno]{amsart}
\usepackage{amsmath,amsthm,amssymb,latexsym,fullpage,setspace,graphicx,float,xcolor,hyperref,verbatim,mathtools}
\usepackage{tikz-cd,tikz,pgfplots}
\usepackage[mode=buildnew]{standalone}
\usetikzlibrary{decorations.markings,math,positioning,graphs}
\usepackage[utf8]{inputenc}
\usepackage[english]{babel}

\numberwithin{equation}{section}

\theoremstyle{plain}
\newtheorem{thm}{Theorem}[section]
\newtheorem*{thm4.8}{Theorem 4.8 of \cite{DKS2025}}
\newtheorem{cor}[thm]{Corollary}
\newtheorem{lem}[thm]{Lemma}

\newtheorem{prop}[thm]{Proposition}

\theoremstyle{definition}
\newtheorem{defn}[thm]{Definition}
\newtheorem{eg}[thm]{Example}

\theoremstyle{remark}
\newtheorem{remark}[thm]{Remark}

\makeatletter
\@namedef{subjclassname@2020}{%
  \textup{2020} Mathematics Subject Classification}
\makeatother

\title{A unifying theory for metrical results on regular continued fraction convergents and mediants}
\author{Karma Dajani}
\address{Department of Mathematics, Utrecht University, P.O.~Box 80010, 3508 TA Utrecht, The Netherlands}
\email{k.dajani@uu.nl}
\author{Cor Kraaikamp}
\address{Delft University of Technology, EWI (DIAM), Mekelweg 4, 2628 CD Delft, The
Netherlands}
\email{c.kraaikamp@tudelft.nl}
\author{Slade Sanderson}
\address{Universit\'e Paris Cit\'e, CNRS, IRIF, F-75013, Paris, France}
\email{slade.sanderson@irif.fr}
\thanks{First published in \emph{Math.~Comp.}~94 (2025), published by the American Mathematical Society.  \textcopyright~2025 American Mathematical Society.}
\date{\today}

\def\convmat#1{
\begin{pmatrix}
u_{#1} & t_{#1}\\
s_{#1} & r_{#1}
\end{pmatrix}
}

\def\convmatR#1{
\begin{pmatrix}
u_{#1}^R & t_{#1}^R\\
s_{#1}^R & r_{#1}^R
\end{pmatrix}
}

\def\F{
\mathcal{F}
}

\begin{document}
\begin{abstract}
We revisit Ito's (\cite{I1989}) natural extension of the Farey tent map, which generates all regular continued fraction convergents and mediants of a given irrational.  With a slight shift in perspective on the order in which these convergents and mediants arise, this natural extension is shown to provide an elegant and powerful tool in the metric theory of continued fractions.  A wealth of old and new results---including limiting distributions of approximation coefficients, analogues of a theorem of Legendre and their refinements, and a generalisation of L\'evy's Theorem to subsequences of convergents and mediants---are presented as corollaries within this unifying theory.
\end{abstract}

\maketitle
\tikzset{->-/.style={decoration={markings,mark=at position #1 with {\arrow{>}}},postaction={decorate}}}

\section{Introduction}\label{Introduction}

Ever since the pioneering work by Shunji Ito, Hitoshi Nakada and Shigeru Tanaka in the late 1970s and early 1980s on the natural extensions of two different classes of continued fraction algorithms (see~\cite{NIT77,N1981}), natural extensions have played a pivotal role in understanding the metric and arithmetic properties of various families of continued fraction algorithms; see for example\footnote{Note this list is far from complete!} Nakada's $\alpha$-expansions (\cite{K1991,LM2008,dJK2018}), the Tanaka-Ito $\alpha$-expansions (\cite{CLS2021,NS2021}), the Katok-Ugarcovici $(a,b)$-continued fractions (\cite{KU2010B,KU2010A,KU2012,AKU2023}), Rosen fractions (\cite{BKS2000,KSS2007,KNS2009,KSS2010}) and its recent generalisations (\cite{CKS2020}). But also for other number theoretic expansions, such as  L\"uroth expansions and $\beta$-expansions (\cite{BBDK1996,DK1996,DKS1996,DK2002A}), the natural extension turned out to be a very important tool. For example, it gives us an alternative method to find the absolutely continuous invariant measure, to prove weak-Bernoullicity of the $\beta$-transformation, and to relate these expansions to L\"uroth-series via inducing (\cite{DKS1996}).

In 1989, Shunji Ito published his groundbreaking paper on mediant convergents (\cite{I1989}). If $x$ is a real irrational number, with \emph{regular continued fraction} ({\sc rcf}) expansion $x=[a_0;a_1,a_2,\dots]$, and with regular convergents $p_n/q_n=[a_0;a_1,\dots,a_n]$ for $n\geq 1$, then the \emph{mediant convergents} of $x$ are given by:
\[
\frac{\lambda p_n+p_{n-1}}{\lambda q_n+q_{n-1}}, \quad \text{for $\lambda\in\mathbb{N}$, $1\leq \lambda < a_{n+1}$}.
\]
\emph{Best approximants}\footnote{This definition of a best approximant is sometimes called a best approximant \emph{of the first kind}.  Best approximants \emph{of the second kind} are defined analogously, but with inequality (\ref{best_approx_ineq}) replaced by $|ex-d|\le |qx-p|$.  Best approximants of the second kind are classified as {\sc rcf}-convergents of $x$; see \cite{RS92}.} of an irrational number $x$ are fractions $p/q$, $q>0,$ in lowest terms, with the property that any other fraction $d/e$  in lowest terms with $e>0$ satisfies
\begin{equation}\label{best_approx_ineq}
\left| x - \frac{d}{e}\right| \leq \left| x - \frac{p}{q}\right| \quad \text{implies that}\quad e>q.
\end{equation}
Now every best approximant is either a regular convergent or a mediant; see~\cite{B1992}.
 
Classically, regular continued fraction convergents of an irrational number $x$ have `excellent' approximation properties. For example, a result by Legendre states that if $p,q\in\mathbb{Z}$ with $q>0,\ \text{gcd}\{ p,q\} =1$ and
\[
\left| x - \frac{p}{q}\right| < \frac{1}{2q^2},
\]
then $p/q$ is a regular continued fraction convergent of $x$, i.e., there exists an $n$ such that $p=p_n$ and $q=q_n$. Furthermore, one can show that if $\theta_n(x):=q_n^2\left| x-\frac{p_n}{q_n}\right|$ for any $n\in\mathbb{N}$, then for every irrational $x$ and every $n\in\mathbb{N}$ one has
\[
\min \{ \theta_{n-1}(x),\theta_n(x),\theta_{n+1}(x)\} < \frac{1}{\sqrt{a_{n+1}^2+4}}\leq \frac{1}{\sqrt{5}}.
\]
For proofs of these results, see~\cite{DK2002B} and the references therein.

In his 1989 paper, Ito studied the natural extension of the so-called \emph{Farey tent map}, an algorithm `underlying' the regular continued fraction expansion which yields all the convergents and mediants of the {\sc rcf}. Ito obtained various metric results on these mediant convergents. In spite of its groundbreaking nature, Ito's paper has generated relatively little attention; see for example~\cite{B1990,J91,BY1996,DK2000,IK2008}. In this paper we want to `repair' this by exploring the possibilities this natural extension of the Farey tent map yields.  In particular, we exploit the fact that the quantities
\[(\lambda q_n+q_{n-1})^2\left|x-\frac{\lambda p_n+p_{n-1}}{\lambda q_n+q_{n-1}}\right|,\quad 0\le \lambda<a_{n+1},\]
may be written explicitly in terms of the forward orbit of $(x,1)$ under Ito's natural extension (see Proposition \ref{Theta&h}).  This fact was essentially known and used in \cite{I1989} and \cite{BY1996}, where the orbit is put into one-to-one correspondence with the sequence of {\sc rcf}-convergents and mediants of $x$, ordered with increasing denominators.  However, we consider a natural rearrangement of this sequence of convergents and mediants which lends itself to a geometrically intuitive one-to-one correspondence with the aforementioned orbit.  This will lead to unified and simple proofs of results from the just mentioned papers~\cite{I1989,B1990,J91,BY1996}, old and classical results by Legendre and Koksma, and various new results such as generalizations of L\'evy's constant and of the Doeblin--Lenstra conjecture to subsequences of convergents and mediants.

This paper is organised as follows: In \S\ref{Background, definitions and notation} we set notation and recall basic facts regarding (semi-)regular continued fractions (\S\ref{(Semi-)regular continued fractions}), the Gauss map underlying {\sc rcf}-expansions (\S\ref{The Gauss map}), the Farey tent map (\S\ref{The Farey tent map}) and the Lehner map and its associated semi-regular continued fraction expansions (\S\ref{Lehner expansions}).  In \S\ref{Farey expansions and convergents}, the Farey tent map is shown to generate semi-regular continued fractions (\S\ref{Farey expansions}) whose convergents (\S\ref{Farey convergents}) consist of all {\sc rcf}-convergents and mediants of a given irrational.  Special emphasis is placed on the order in which these \emph{Farey convergents} arise, and which \emph{differs} from that historically studied in the literature.  In an appendix (\S\ref{Appendix}), it is shown that these two orderings do not affect the statements of central results in \S\ref{Metrical results}, and thus a number of old results are re-obtained from our approach.  Section \ref{Ito's natural extension of the Farey tent map} recalls Ito's natural extension of the Farey tent map and the relationship between its dynamics and Farey convergents.  Subsection \ref{Inducing Ito's natural extension} sets the framework for inducing Ito's natural extension to obtain desired subsequences of Farey convergents, and in \S\ref{The (relative) equidistribution of orbits} we prove central results on the equidistribution of the orbit of $(x,1)$ under the induced maps.  The main metrical results---old and new---are contained in \S\ref{Metrical results}.  Here we consider limiting distributions of approximation coefficients (\S\ref{Approximation coefficients and their limiting distributions}), Legendre-type theorems (\S\ref{On the theorems of Legendre, Fatou--Grace and Koksma}), consecutive approximation coefficients (\S\ref{Consecutive approximation coefficients}) and generalisations of results of L\'evy (\S\ref{A generalised Levy-type theorem}).

\bigskip

{\flushleft \textbf{Acknowledgments.}}  We thank the anonymous referees whose comments greatly improved the exposition of this paper.  This work is part of project number 613.009.135 of the research programme Mathematics Clusters which is financed by the Dutch Research Council (NWO).

\bigskip
\section{Background, definitions and notation}\label{Background, definitions and notation}

\subsection{(Semi-)regular continued fractions}\label{(Semi-)regular continued fractions}
A \textit{semi-regular continued fraction} ({\sc srcf}) is a formal (infinite or finite) expression of the form
\[[\beta_0;\alpha_1/\beta_1,\alpha_2/\beta_2,\dots]=\beta_0+\cfrac{\alpha_1}{\beta_1+\cfrac{\alpha_2}{\beta_2+\ddots}}\]
with $\beta_0\in\mathbb{Z}$, and for each $n\ge 1,\ \alpha_n=\pm 1$ and $\beta_n\ge 1$ integers satisfying 
\[\alpha_{n+1}+\beta_n\ge1,\]
and---in the infinite case---
\[\alpha_{n+1}+\beta_n\ge2\]
infinitely often.  Set
\[B_0=B_0([\beta_0;\alpha_1/\beta_1,\dots]):=\begin{pmatrix}1 & \beta_0\\ 0 & 1\end{pmatrix}\quad \text{and}\quad B_n=B_n([\beta_0;\alpha_1/\beta_1,\dots]):=\begin{pmatrix}0 & \alpha_n\\ 1 & \beta_n\end{pmatrix},\quad n>0,\]
and for $0\le i\le j$,
\[B_{[i,j]}=B_{[i,j]}([\beta_0;\alpha_1/\beta_1,\dots]):=B_iB_{i+1}\cdots B_j.\]
For a matrix $A=\left(\begin{smallmatrix}a & b\\ c & d\end{smallmatrix}\right)$, we denote by $A\cdot z:=\frac{az+b}{cz+d},\ z\in\mathbb{R}\cup\{\infty\},$ the action of $A$ as a M\"obius transformation.  Writing the entries of $B_{[0,n]}$ as $\left(\begin{smallmatrix}R_n & P_n\\ S_n & Q_n\end{smallmatrix}\right)$, the fraction 
\[\frac{P_n}{Q_n}:=B_{[0,n]}\cdot 0=\beta_0+\cfrac{\alpha_1}{\beta_1+\cfrac{\alpha_2}{\ddots+\cfrac{\alpha_n}{\beta_n}}}=[\beta_0;\alpha_1/\beta_1,\dots,\alpha_n/\beta_n]\in\mathbb{Q}\]
is called the \textit{$n^\text{th}$ convergent} of $[\beta_0;\alpha_1/\beta_1,\alpha_2/\beta_2,\dots]$.  By Tietze's Convergence Theorem (see, say, \cite{perron_54}) the above conditions on the digits $\alpha_n$ and $\beta_n$ guarantee that $x=\lim_{n\to\infty}\frac{P_n}{Q_n}\in\mathbb{R}$ always exists, and thus we call $[\beta_0;\alpha_1/\beta_1,\alpha_2/\beta_2,\dots]$ a \textit{{\sc srcf}-expansion} of $x$ and refer to the convergents $P_n/Q_n$ as convergents of $x$.\footnote{We emphasise that a real number $x$ has many {\sc srcf}-expansions, and the convergents of $x$ depend on the expansion in question.}  The convergents $P_n/Q_n$ of any {\sc srcf}-expansion of $x\in\mathbb{R}$ are reduced, as $\det (B_0)=1$ and $\det (B_{[0,n]}) = \alpha_1\cdots\alpha_n(-1)^n\in \{\pm 1\}$, for $n\geq 1$.  Notice for any $n\ge 0$ that
\[\begin{pmatrix}R_{n+1} & P_{n+1}\\ S_{n+1} & Q_{n+1}\end{pmatrix}=B_{[0,n]}B_{n+1}=\begin{pmatrix}R_{n} & P_{n}\\ S_{n} & Q_{n}\end{pmatrix}\begin{pmatrix}0 & \alpha_{n+1}\\ 1 & \beta_{n+1}\end{pmatrix}=\begin{pmatrix}P_n & \beta_{n+1}P_n+\alpha_{n+1}R_n\\ Q_n & \beta_{n+1}Q_n+\alpha_{n+1}S_n\end{pmatrix}.\]
In particular, $R_{n+1}=P_n$ and $S_{n+1}=Q_n$.  In view of the definition of $B_{[0,0]}=B_0$, we set $P_{-1}:=R_0=1$ and $Q_{-1}:=S_0=0$, and call $P_{-1}/Q_{-1}=1/0=\infty$ the \textit{$(-1)^\text{st}$ convergent} of $[\beta_0;\alpha_1/\beta_1,\alpha_2/\beta_2,\dots]$ and of $x$.  This gives the following recurrence relations for all $n\ge 0$:
\begin{alignat}{3}\label{rec_rels}
P_{n+1}&=\beta_{n+1}P_n+\alpha_{n+1}P_{n-1}, \qquad &P_{-1}=1,\ &P_0=\beta_0,\\
Q_{n+1}&=\beta_{n+1}Q_n+\alpha_{n+1}Q_{n-1}, \qquad &Q_{-1}=0,\ &Q_0=1.\notag
\end{alignat}
 
A \textit{regular continued fraction} ({\sc rcf}) is a semi-regular continued fraction with $\alpha_n=1,\ n\ge 1$ (note that the conditions of a {\sc srcf} are now trivially satisfied for any $\beta_n\ge 1$).  A {\sc rcf} is also denoted by
\[[a_0;a_1,a_2,\dots]:=[a_0;1/a_1,1/a_2,\dots]\]
and its sequence of convergents by $(p_n/q_n)_{n\ge -1}$.  The digit $a_n$ is the \emph{$n^\text{th}$ partial quotient} of $x=[a_0;a_1,a_2,\dots]$.  The \textit{mediant convergents} (or, simply, \emph{mediants}) of $x$ are defined as the fractions 
\[\frac{\lambda p_n+p_{n-1}}{\lambda q_n+q_{n-1}},\quad \text{for $\lambda\in\mathbb{N}$, $1\leq \lambda < a_{n+1}$}.\]

\subsection{The Gauss map}\label{The Gauss map}
The \textit{Gauss map} $G:[0,1]\to [0,1]$ is defined by $G(0)=0$ and $G(x)=1/x-\lfloor1/x\rfloor,\ x>0$ (see Figure \ref{GaussFarey}).\footnote{While $G$ may also be defined as a self-map of $[0,1)$, we choose to include the endpoint $1$ in our definition for later notational purposes.}  For $x\in\mathbb{R}$, let $a_0=a_0(x):=\lfloor x\rfloor$ and $x_0:=x-a_0\in[0,1)$.  Define $a(x):=\lfloor 1/x\rfloor,\ x\neq 0$, and set $a_n=a_n(x):=a(G^{n-1}(x_0)),\ G^{n-1}(x_0)\neq 0$.  With this notation, for $G^{n-1}(x_0)\neq 0$,
\[G^n(x_0)=\frac1{G^{n-1}(x_0)}-a_n,\]
which can be rewritten as
\[G^{n-1}(x_0)=\frac{1}{a_n+G^n(x_0)}.\] 
Repeatedly applying this last relation, one finds that
\begin{equation}\label{expn_eqn}
x=a_0+\cfrac{1}{a_1+\cfrac{1}{a_2+\ddots+\cfrac{1}{a_n+G^n(x_0)}}}=[a_0;a_1,\dots,a_{n-1},a_n+G^n(x_0)].
\end{equation}
From the Euclidean algorithm it follows that for every rational $x$ there exists an $n\geq 0$ such that $G^n(x_0) = 0$, and we see that the {\sc rcf}-expansion of $x\in\mathbb{Q}$ generated by $G$ is finite.  In fact, every rational $x$ has precisely two {\sc rcf}-expansions, the first of which is generated by $G$:
\begin{equation}\label{two_expns}
x=[a_0;a_1,\dots,a_n]\quad \text{and}\quad x=[a_0;a_1,\dots,a_n-1,1],
\end{equation}
where $a_n\ge 2$ when $n\ge 1$ (that is, when $x\notin\mathbb{Z}$).  The \emph{depth} of $x\in\mathbb{Q}$ is the number $n$ occurring in the above two expansions.  

For $x\in\mathbb{R}\backslash\mathbb{Q}$, taking $n\to\infty$ in (\ref{expn_eqn}), we see that the Gauss map $G$ generates a (unique) {\sc rcf}-expansion $[a_0;a_1,a_2,\dots]$ of $x$.  It is well known that the dynamical system $([0,1],\mathcal{B},\nu_G,G)$ is ergodic, where $\mathcal{B}$ denotes the Borel $\sigma$-algebra and $\nu_G$ is the \textit{Gauss measure} with density $\frac1{\log 2 (1+x)}$ (see, say, \cite{DK2002B}).

Let $\Omega:=[0,1]^2$, and define $\mathcal{G}:\Omega\to \Omega$ by $\mathcal{G}(0,y)=(0,y)$ and, for $x\neq 0$,
\[\mathcal{G}(x,y):=\left(G(x),\frac1{a(x)+y}\right).\]
For $(x,y)\in\Omega,\ x\in(0,1),$ with {\sc rcf}-expansions $(x,y)=([0;a_1,a_2,\dots],[0;b_1,b_2,\dots])$, the map $\mathcal{G}$ acts as a two-dimensional shift 
\begin{equation}\label{GaussNE}
\mathcal{G}(x,y)=([0;a_2,a_3,a_4,\dots],[0;a_1,b_1,b_2,\dots]).
\end{equation}
In \cite{NIT77,N1981}, the authors show that the ergodic system $(\Omega,\mathcal{B},\bar\nu_G,\mathcal{G})$ is the natural extension of $([0,1],\mathcal{B},\nu_G,G)$, where $d\bar\nu_G=\frac{dxdy}{\log 2(1+xy)^2}$ (see also \cite{DK2002B,dajani_kalle_21}).    In fact, $(\Omega,\mathcal{B},\bar\nu_G,\mathcal{G})$ and $([0,1],\mathcal{B},\nu_G,G)$ are strongly mixing (\cite{IK02}).

\subsection{The Farey tent map}\label{The Farey tent map}

\begin{figure}[t]
\includestandalone[width=.4\textwidth]{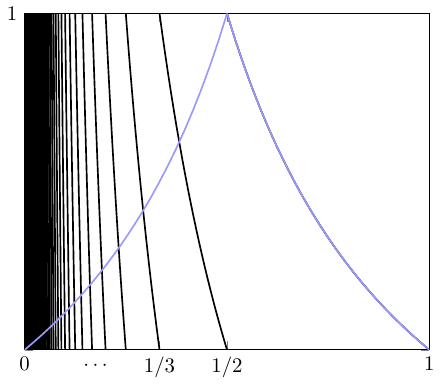}
\caption{The Gauss map $G$ (black) and the Farey tent map $F$ (blue).  Both maps coincide on the domain $[1/2,1]$.}
\label{GaussFarey}
\end{figure}

Let $F:[0,1]\to[0,1]$ denote the \textit{Farey tent map} given by
\[
F(x):=\begin{cases}
x/(1-x), & x\le1/2,\\
(1-x)/x, & x> 1/2;
\end{cases}
\]
see Figure \ref{GaussFarey}.  The dynamical system $([0,1],\mathcal{B},\mu,F)$ is ergodic, where $\mu$ is the absolutely continuous, infinite, $\sigma$-finite, $F$-invariant measure with density $1/x$ (\cite{daniels_62,parry_62,I1989}).  It follows from the definition of $F$ that if $x\in(0,1)$ has {\sc rcf}-expansion $x=[0;a_1,a_2,a_3\dots]$, then
\begin{equation}\label{symbolic_farey}
F(x)=\begin{cases}
[0;a_1-1,a_2,a_3,\dots], & a_1>1,\\
[0;a_2,a_3,a_4,\dots], & a_1=1.
\end{cases}
\end{equation}
From this, one finds that the Gauss map $G$ is the jump transformation of $F$ associated to the interval $(1/2,1]$; in particular, for $x$ as above,
\[\min\{j\ge0\ |\ F^j(x)\in(1/2,1]\}=a_1-1, \quad \text{and} \quad G(x)=F^{a_1}(x).\]

Define $\varepsilon:[0,1]\to\{0,1\}$ by
\[
\varepsilon(x):=\begin{cases}
0, & x\le1/2,\\
1, & x> 1/2,
\end{cases}
\]
and for $x\in[0,1]$ and $n\ge 1$, set $\varepsilon_n=\varepsilon_n(x):=\varepsilon(F^{n-1}(x))$.  From (\ref{symbolic_farey}), it follows that for $x=[0;a_1,a_2,\dots]$,
\begin{equation}\label{epsilon_seq}
\varepsilon_1\varepsilon_2\varepsilon_3\cdots=0^{a_1-1}10^{a_2-1}10^{a_3-1}1\cdots
\end{equation}
(see also \cite{I1989}).

\subsection{Lehner expansions}\label{Lehner expansions}

It was shown in \cite{I1989} that the Farey tent map generates all convergents and mediant convergents of the {\sc rcf}-expansion of any irrational $x\in (0,1)$. Originally there was no continued fraction algorithm `attached' to the Farey tent map $F$. Such a continued fraction expansion does exist, and can be obtain via the so-called Lehner map $L$, which was introduced by Joe Lehner in 1994; see \cite{lehner_94}, and also \cite{DK2000}.

Let $L:[1,2]\to[1,2]$ be given by
\[
L(x):=\begin{cases}
1/(2-x), & x\le3/2,\\
1/(x-1), & x>3/2,
\end{cases}
\]
and for $x\in[1,2]$ and each $n\ge 0$, set
\[
(b_n,e_{n+1})=(b_n(x),e_{n+1}(x)):=\begin{cases}
(2,-1), & L^n(x)\le3/2,\\
(1,1), & L^n(x)> 3/2.
\end{cases}
\]
The digits $(b_n,e_{n+1})$ generate the so-called \textit{Lehner expansion} of $x\in[1,2]$,
\begin{equation}\label{lehner_expn}
x=[b_0;e_1/b_1,e_2/b_2,\dots],
\end{equation}
which is a {\sc srcf}-expansion (see \cite{lehner_94,DK2000}).

\bigskip
\section{Farey expansions and convergents}\label{Farey expansions and convergents}

\subsection{Farey expansions}\label{Farey expansions}

As noted above, Ito (\cite{I1989}) studied the ergodic properties of the dynamical system $([0,1],\mathcal{B},\mu,F)$ without any explicit mention of associated (semi-regular) continued fraction expansions, while Lehner (\cite{lehner_94}) studied expansions of the form (\ref{lehner_expn}) generated by $L$ but no dynamical properties of this map.  In \cite{DK2000} it is observed that the dynamical systems $([0,1],\mathcal{B},\mu,F)$ and $([1,2],\mathcal{B},\rho,L)$ are isomorphic via the translation $x\mapsto x+1$, where $d\rho=dx/(x-1)$.  Via this isomorphism, the map $F$ can be used to generate a \textit{Farey expansion} for each $x\in[0,1]$ (see also \cite{iosifescu_sebe_08}).  Indeed, fix $x\in[0,1]$, and let $[b_0;e_1/b_1,e_2/b_2,\dots]$ be the Lehner expansion of $x+1$.  Then $[b_0-1;e_1/b_1,e_2/b_2,\dots]$ is a {\sc srcf}-expansion of $x$, and we find that
\begin{align*}
(b_n,e_{n+1})&=\left.\begin{cases}
(2,-1), & L^n(x+1)\le3/2\\
(1,1), & L^n(x+1)> 3/2
\end{cases}\right\}
=\left.\begin{cases}
(2,-1), & F^n(x)\le1/2\\
(1,1), & F^n(x)> 1/2
\end{cases}\right\}\\
&=\left.\begin{cases}
(2,-1), & \varepsilon_{n+1}=0\\
(1,1), & \varepsilon_{n+1}=1
\end{cases}\right\}
=(2-\varepsilon_{n+1},2\varepsilon_{n+1}-1).
\end{align*}
Hence $F$ generates the \textit{Farey expansion} 
\[x=[b_0-1;e_1/b_1,e_2/b_2,\dots],\] 
where
\begin{equation}\label{farey_digits}
(b_n,e_{n+1})=(2-\varepsilon_{n+1},2\varepsilon_{n+1}-1),\quad n\ge 0.
\end{equation}
The corresponding convergents $P_n/Q_n=[b_0-1;e_1/b_1,e_2/b_2,\dots,e_n/b_n]$ are called the \textit{Farey convergents} of $x$.

\subsection{Farey convergents}\label{Farey convergents}

Fix an irrational $x\in(0,1)$ with {\sc rcf}-expansion $x=[0;a_1,a_2,\dots]$ and convergents $p_n/q_n$.  In \cite{I1989} and \cite{DK2000} it is shown that $F$ (respectively $L$) generates all convergents and mediants of the {\sc rcf}-expansion of $x$ (respectively $x+1$).  We reproduce this fact for $F$ here, fixing notation\footnote{Notation is largely recycled from \cite{I1989}; however, matrix entries are permuted so as to conform with modern notation of their action via M\"obius transformations.} along the way and paying special attention to the order in which these convergents and mediants arise.  

Set
\[A_0:=\begin{pmatrix}1 & 0\\ 1 & 1\end{pmatrix}\quad \text{and}\quad A_1:=\begin{pmatrix}0 & 1\\ 1 & 1\end{pmatrix},\]
or, more succinctly, 
\begin{equation}\label{A_epsilon}
A_\varepsilon:=\begin{pmatrix}1-\varepsilon & \varepsilon\\ 1 & 1\end{pmatrix},\quad \varepsilon\in\{0,1\}.
\end{equation}
Note that as M\"obius transformations, $A_0^{-1}$ and $A_1^{-1}$ correspond to the left and right branches of $F$, respectively.  In particular, $F(x)=A_{\varepsilon(x)}^{-1}\cdot x,$ so $x=A_{\varepsilon(x)}\cdot F(x)$.  Setting $x_n:=F^n(x)$ for $n\ge 0$, we find that 
\[x=(A_{\varepsilon_1}A_{\varepsilon_2}\cdots A_{\varepsilon_{n}})\cdot x_n,\]
where $\varepsilon_n=\varepsilon_n(x)=\varepsilon(F^{n-1}(x))$ for all $n\ge 1$ (see \S\ref{The Farey tent map}).  For $0\le i\le j$, define
\begin{equation}\label{A_[0,n]}
A_{[i,j]}=A_{[i,j]}(x):=A_{\varepsilon_i(x)}A_{\varepsilon_{i+1}(x)}\cdots A_{\varepsilon_j(x)} \quad \text{and} \quad A_{[j,i]}=A_{[j,i]}(x):=A_{\varepsilon_j(x)}A_{\varepsilon_{j-1}(x)}\cdots A_{\varepsilon_i(x)},
\end{equation}
where $A_{\varepsilon_0(x)}:=I_2$ is the two-by-two identity matrix.  Denote the entries of $A_{[0,n]},\ n\ge 0,$ by
\[
\convmat{n}=
\begin{pmatrix}
u_n(x) & t_n(x)\\
s_n(x) & r_n(x)
\end{pmatrix}
:=A_{[0,n]}.
\]
Observe that for any $k\in\mathbb{Z}$,
\begin{equation}\label{A_0^kA_1}
A_0^kA_1=\begin{pmatrix}1 & 0\\ 1 & 1\end{pmatrix}^k\begin{pmatrix}0 & 1\\ 1 & 1\end{pmatrix}=\begin{pmatrix}1 & 0\\ k & 1\end{pmatrix}\begin{pmatrix}0 & 1\\ 1 & 1\end{pmatrix}=\begin{pmatrix}0 & 1\\ 1 & k+1\end{pmatrix}.
\end{equation}
In view of (\ref{epsilon_seq}), for $n\ge 0$ we set
\begin{equation}\label{j_lambda}
j_n=j_n(x):=\#\{1\le k\le n\ |\ \varepsilon_k=1\}\quad \text{and}\quad \lambda_n=\lambda_n(x):=n-\sum_{k=1}^{j_n} a_k.
\end{equation}
That is, $j_n$ and $\lambda_n$ are the unique integers satisfying
\begin{equation}\label{n_exp_with_x_digits}
n=a_1+a_2+\dots+a_{j_n}+\lambda_n,\qquad j_n\ge 0,\quad 0\le \lambda_n<a_{j_n+1}.
\end{equation}
From (\ref{epsilon_seq}) and (\ref{A_0^kA_1}), it follows for $n>0$ that
\begin{align}\label{convmats}
\convmat{n}=A_{[0,n]}&=I_2A_{\varepsilon_1}\cdots A_{\varepsilon_n}\notag \\
&=A_0^{a_1-1}A_1\cdots A_0^{a_{j_n}-1}A_1A_0^{\lambda_n} \notag\\
&=\begin{pmatrix}0 & 1\\ 1 & a_1\end{pmatrix}\cdots\begin{pmatrix}0 & 1\\ 1 & a_{j_n}\end{pmatrix}\begin{pmatrix}1 & 0\\ {\lambda_n} & 1\end{pmatrix} \notag\\
&=\begin{pmatrix}p_{{j_n}-1} & p_{j_n}\\ q_{{j_n}-1} & q_{j_n}\end{pmatrix}\begin{pmatrix}1 & 0\\ {\lambda_n} & 1\end{pmatrix} \notag\\
&=\begin{pmatrix}{\lambda_n} p_{j_n}+p_{{j_n}-1} & p_{j_n}\\ {\lambda_n} q_{j_n}+q_{{j_n}-1} & q_{j_n}\end{pmatrix}
\end{align}
(see also Lemma 1.1 of \cite{I1989}).  Note that equality of the first and final expressions also holds for $n={j_n}={\lambda_n}=0$, for in this case both matrices are the identity $I_2$.  We also have for $n>0$ that
\begin{align}\label{convmats_reversed}
A_{[n,0]}&=A_{\varepsilon_n}\cdots A_{\varepsilon_1}I_2\notag\\
&=A_0^{\lambda_n} A_1A_0^{a_{j_n}-1}\cdots A_1A_0^{a_1-1}A_1A_1^{-1}\notag\\
&=\begin{pmatrix}0 & 1\\ 1 & {\lambda_n}+1\end{pmatrix}\begin{pmatrix}0 & 1\\ 1 & a_{j_n}\end{pmatrix}\cdots\begin{pmatrix}0 & 1\\ 1 & a_1\end{pmatrix}\begin{pmatrix}-1 & 1\\ 1 & 0\end{pmatrix}\notag\\
&=\begin{pmatrix}0 & 1\\ 1 & {\lambda_n}+1\end{pmatrix}\left(\begin{pmatrix}0 & 1\\ 1 & a_1\end{pmatrix}\cdots\begin{pmatrix}0 & 1\\ 1 & a_{j_n}\end{pmatrix}\right)^T\begin{pmatrix}-1 & 1\\ 1 & 0\end{pmatrix}\notag\\
&=\begin{pmatrix}0 & 1\\ 1 & {\lambda_n}+1\end{pmatrix}\begin{pmatrix}p_{j_n-1} & q_{j_n-1}\\ p_{j_n} & q_{j_n}\end{pmatrix}\begin{pmatrix}-1 & 1\\ 1 & 0\end{pmatrix}\notag\\
&=\begin{pmatrix}q_{j_n}-p_{j_n} & p_{j_n}\\ ({\lambda_n}+1)q_{j_n}+q_{j_n-1}-(({\lambda_n}+1)p_{j_n}+p_{j_n-1}) &({\lambda_n}+1)p_{j_n}+p_{j_n-1}\end{pmatrix}\notag\\
&=\begin{pmatrix}r_n-t_n & t_n\\ s_n+r_n-(u_n+t_n) & u_n+t_n\end{pmatrix}
\end{align}
(cf. (\ref{convmats})), and the first and final expressions are again also equal to $I_2$ for $n=0$.
From (\ref{convmats}) it is clear that the set $\{u_n/s_n\}_{n\ge 0}$ equals the set 
\[\left\{\frac{\lambda p_j+p_{j-1}}{\lambda q_j+q_{j-1}}\ \Big|\ j\ge 0,\ 0\le\lambda<a_{j+1}\right\}\]
of all convergents and mediants of the {\sc rcf}-expansion $[0;a_1,a_2,\dots]$ of $x$.  In fact, the sequence $(u_n/s_n)_{n\ge 0}$ is precisely the sequence of Farey convergents of $x$:

\begin{prop}\label{FareyConvs}
For each $n\ge 0$,
\[\begin{pmatrix}u_n \\ s_n\end{pmatrix}=\begin{pmatrix}P_{n-1} \\ Q_{n-1}\end{pmatrix},\]
where $P_n/Q_n=[b_0-1;e_1/b_1,\dots,e_n/b_n]$ is the $n^{\text{th}}$ Farey convergent of $x$.
\end{prop}
\begin{proof}
The proof is by induction.  For $n=0$, we have
\[\begin{pmatrix}u_0 \\ s_0\end{pmatrix}=\begin{pmatrix}1 \\ 0\end{pmatrix}=\begin{pmatrix}P_{-1} \\ Q_{-1}\end{pmatrix},\]
and for $n=1$,
\[\begin{pmatrix}u_1 \\ s_1\end{pmatrix}=\begin{pmatrix}1-\varepsilon_1 \\ 1\end{pmatrix}=\begin{pmatrix}(2-\varepsilon_1)-1 \\ 1\end{pmatrix}=\begin{pmatrix}b_0-1 \\ 1\end{pmatrix}=\begin{pmatrix}P_{0} \\ Q_{0}\end{pmatrix}.\]
Now suppose for some $n\ge 1$ it holds for each $k\le n$ that
\[\begin{pmatrix}u_k \\ s_k\end{pmatrix}=\begin{pmatrix}P_{k-1} \\ Q_{k-1}\end{pmatrix}.\]
Now
\begin{align}\label{farey_conv_mat}
\convmat{n}=\convmat{n-1}\begin{pmatrix}1-\varepsilon_n & \varepsilon_n \\ 1 & 1\end{pmatrix}=\begin{pmatrix}(1-\varepsilon_n)u_{n-1}+t_{n-1} & \varepsilon_nu_{n-1}+t_{n-1} \\ (1-\varepsilon_n)s_{n-1}+r_{n-1} & \varepsilon_ns_{n-1}+r_{n-1}\end{pmatrix}.
\end{align}
The top-left entries give
\[t_{n-1}=u_n+(\varepsilon_n-1)u_{n-1}.\]
Then the top-right entries---replacing $t_{n-1}$ with the right-hand side of the previous line---give
\[t_n=u_n+(2\varepsilon_n-1)u_{n-1}.\]
Since (\ref{farey_conv_mat}) holds for all $n\ge 1$, the top-left entries (replacing $n$ with $n+1$) together with the previous line, (\ref{farey_digits}) and (\ref{rec_rels}) gives
\[u_{n+1}=(1-\varepsilon_{n+1})u_{n}+t_{n}=(2-\varepsilon_{n+1})u_{n}+(2\varepsilon_n-1)u_{n-1}=b_nP_{n-1}+e_nP_{n-2}=P_n.\]
With similar computations one finds $s_{n+1}=Q_n$.
\end{proof}

As a sequence, (\ref{convmats}) gives 
\begin{alignat*}{3}
\left(\convmat{n}\right)_{n\ge 0}
=\bigg(&\begin{pmatrix}p_{-1} & p_0 \\ q_{-1} & q_0\end{pmatrix}&&,\begin{pmatrix}p_0+p_{-1} & p_0 \\ q_0+q_{-1} & q_0\end{pmatrix}&&,\dots,\begin{pmatrix}(a_1-1)p_0+p_{-1} & p_0 \\ (a_1-1)q_0+q_{-1} & q_0\end{pmatrix},\\
&\begin{pmatrix}p_{0} & p_1 \\ q_{0} & q_1\end{pmatrix}&&,\begin{pmatrix}p_1+p_{0} & p_1 \\ q_1+q_0 & q_1\end{pmatrix}&&,\dots,\begin{pmatrix}(a_2-1)p_1+p_0 & p_1 \\ (a_2-1)q_1+q_0 & q_1\end{pmatrix},\dots,\\
&\begin{pmatrix}p_{j-1} & p_j \\ q_{j-1} & q_j\end{pmatrix}&&,\begin{pmatrix}p_j+p_{j-1} & p_j \\ q_j+q_{j-1} & q_j\end{pmatrix}&&,\dots,\begin{pmatrix}(a_{j+1}-1)p_j+p_{j-1} & p_j \\ (a_{j+1}-1)q_j+q_{j-1} & q_j\end{pmatrix},\dots\bigg),
\end{alignat*}
and thus by Proposition \ref{FareyConvs}, the Farey convergents occur in the following order:  
\begin{alignat}{3}\label{FLseq}
\left(\frac{P_{n-1}}{Q_{n-1}}\right)_{n\ge0}=\left(\frac{u_n}{s_n}\right)_{n\ge0}
=\bigg(&\frac{p_{-1}}{q_{-1}}&&,\frac{p_0+p_{-1}}{q_0+q_{-1}}&&,\dots,\frac{(a_1-1)p_0+p_{-1}}{(a_1-1)q_0+q_{-1}},\\
&\frac{p_{0}}{q_{0}}&&,\frac{p_1+p_0}{q_1+q_0}&&,\dots,\frac{(a_2-1)p_1+p_0}{(a_2-1)q_1+q_0},\dots,\notag\\
&\frac{p_{j-1}}{q_{j-1}}&&,\frac{p_j+p_{j-1}}{q_j+q_{j-1}}&&,\dots,\frac{(a_{j+1}-1)p_j+p_{j-1}}{(a_{j+1}-1)q_j+q_{j-1}},\dots\bigg).\notag
\end{alignat}
Notice that the denominators $(s_n)_{n\ge 0}$ do not form an increasing sequence.  Supposedly to `remedy' this, in \cite{I1989} Ito instead considers the collection $\{(u_n+t_n)/(s_n+r_n)\}_{n\ge 0}$.  As a sequence, this gives
\begin{alignat*}{4}
\left(\frac{u_n+t_n}{s_n+r_n}\right)_{n\ge0}
=\bigg(&\frac{p_0+p_{-1}}{q_0+q_{-1}}&&,\frac{2p_0+p_{-1}}{2q_0+q_{-1}}&&,\dots,\frac{(a_1-1)p_0+p_{-1}}{(a_1-1)q_0+q_{-1}}&&,\frac{p_1}{q_1},\\
&\frac{p_1+p_{0}}{q_1+q_{0}}&&,\frac{2p_1+p_0}{2q_1+q_0}&&,\dots,\frac{(a_2-1)p_1+p_0}{(a_2-1)q_1+q_0}&&,\frac{p_2}{q_2},\dots,\\
&\frac{p_j+p_{j-1}}{q_j+q_{j-1}}&&,\frac{2p_j+p_{j-1}}{2q_j+q_{j-1}}&&,\dots,\frac{(a_{j+1}-1)p_j+p_{j-1}}{(a_{j+1}-1)q_j+q_{j-1}}&&,\frac{p_{j+1}}{q_{j+1}},\dots\bigg).
\end{alignat*}
However, in light of Proposition \ref{FareyConvs} and results in the sequel, we find it more natural to study the Farey convergents $u_n/s_n$.

\bigskip
\section{Ito's natural extension of the Farey tent map}\label{Ito's natural extension of the Farey tent map}

\begin{figure}
\includestandalone[width=.2\textwidth]{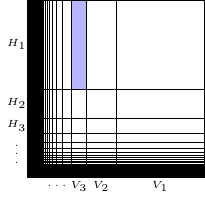}
\includestandalone[width=.2\textwidth]{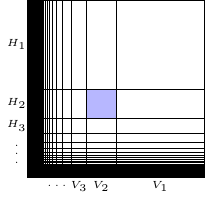}
\includestandalone[width=.2\textwidth]{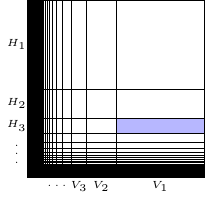}
\includestandalone[width=.2\textwidth]{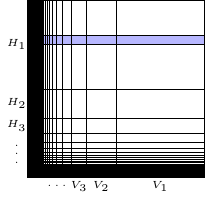}
\caption{From left to right: The sets $V_3\cap H_1,\ \F(V_3\cap H_1),\ \F^2(V_3\cap H_1)$ and $\F^3(V_3\cap H_1)$, respectively.}
\label{FareyNE}
\end{figure}

In \cite{I1989}, Ito determined a planar natural extension $(\Omega,\mathcal{B},\bar\mu,\F)$ of the dynamical system $([0,1],\mathcal{B},\mu,F)$ associated to the Farey tent map.  The map $\F:\Omega\to\Omega$ is given by
\begin{equation}\label{FareyNEeqn}
\F(x,y):=\begin{cases}
\left(\frac{x}{1-x},\frac{y}{1+y}\right), & x\le1/2,\\
\left(\frac{1-x}{x},\frac{1}{1+y}\right), & x> 1/2,
\end{cases}
\end{equation}
where again $\Omega=[0,1]^2$, and $\bar\mu$ is the absolutely continuous measure on $\Omega$ with density $1/(x+y-xy)^2$.  The measure $\bar\mu$ is infinite, $\sigma$-finite and $\F$-invariant, and the natural extension $(\Omega,\mathcal{B},\bar\mu,\F)$ is ergodic (Theorem 1.3 of \cite{I1989}).  Using the matrix notation from (\ref{A_epsilon}), $\F$ may be written as 
\[\F(x,y)=\left(A_{\varepsilon(x)}^{-1}\cdot x,A_{\varepsilon(x)}\cdot y\right),\]
and thus the $n^{\text{th}}$ iterate is
\[\F^n(x,y)=\left(A_{[0,n]}^{-1}\cdot x,A_{[n,0]}\cdot y\right)\]
(recall (\ref{A_[0,n]})).  The map $\F$ admits a particularly nice geometric interpretation, which we now describe.  For each integer $k\ge 1$, let 
\[V_k:=\left(\frac1{k+1},\frac1{k}\right]\times [0,1]\quad \text{and}\quad H_k:=[0,1]\times\left(\frac1{k+1},\frac1{k}\right]\]
denote the \textit{$k^{\textit{th}}$ vertical} and \textit{horizontal regions}, respectively.  Now fix $(x,y)\in\Omega$ with {\sc rcf}-expansions 
\[(x,y)=([0;a_1,a_2,\dots],[0;b_1,b_2,\dots])\in V_{a_1}\cap H_{b_1}.\]
One verifies using (\ref{symbolic_farey}) and (\ref{FareyNEeqn}) that for $x\neq 1$,
\begin{equation}\label{F_NE_symbolic}
\F(x,y)=\begin{cases}
([0;a_1-1,a_2,\dots],[0;b_1+1,b_2,\dots]), & a_1>1,\\
([0;a_2,a_3\dots],[0;1,b_1,b_2,\dots]), & a_1=1.
\end{cases}
\end{equation}
Thus the image of the rectangle $V_{a}\cap H_{b},\ a>1,$ is the rectangle $\F(V_{a}\cap H_{b})=V_{a-1}\cap H_{b+1}$ immediately below and to the right of the original rectangle, and the image of the right-half $V_1$ of $\Omega$ is the top half $\F(V_1)=H_1$, modulo a Lebesgue-null set.  In particular, subsequent iterates $\F^\lambda,\ 0\le \lambda<a$, `slide' the rectangle $V_a\cap H_1$ `diagonally' along $a$ rectangles, and the next iterate $\F^{a}(V_{a}\cap H_1)$ is mapped back as a subset of $H_1$ (see Figure \ref{FareyNE}).  

For $x=[0;a_1,a_2,\dots]\in(0,1)\backslash\mathbb{Q}$ and $n\ge 0$ set 
\begin{equation}\label{(xn,yn)}
(x_n,y_n):=\F^n(x,1)=\left(A_{[0,n]}^{-1}\cdot x,A_{[n,0]}\cdot 1\right).
\end{equation} 
The above geometric interpretation of $\mathcal{F}$ provides a natural identification between the orbit $(x_n,y_n)_{n\ge 0}$ in $\Omega$ and the sequence $(P_{n-1}/Q_{n-1})_{n\ge0}=(u_n/s_n)_{n\ge 0}$ of Farey convergents of $x$ from (\ref{FLseq}).  In particular, the first $a_1$ points
\[(x_0,y_0),(x_1,y_1),\dots,(x_{a_1-1},y_{a_1-1})\]
in the orbit belong to the rectangles 
\[V_{a_1}\cap H_1,V_{a_1-1}\cap H_2,\dots,V_1\cap H_{a_1}\] 
and correspond to the Farey convergents 
\[\frac{p_{-1}}{q_{-1}},\frac{p_0+p_{-1}}{q_0+q_{-1}},\dots,\frac{(a_1-1)p_0+p_{-1}}{(a_1-1)q_0+q_{-1}},\]
respectively.  Note that $(x_{a_1},y_{a_1})=([0;a_2,a_3,\dots],[0;1,a_1])$, so the next $a_2$ points 
\[(x_{a_1},y_{a_1}),(x_{a_1+1},y_{a_1+1}),\dots,(x_{a_1+a_2-1},y_{a_1+a_2-1})\]
of the orbit belong to (the closures\footnote{The closure is needed if and only if $a_1=1$, since then $y_{a_1+\lambda}=[0;\lambda+1,a_1]=1/(\lambda+2)\notin H_{\lambda+1}$ for $0\le \lambda<a_2$.  As shown below, this annoyance is `corrected' for $n\ge a_1+a_2$, and the closures are no longer needed.  Throughout the paper, we shall overlook this innocuous subtlety and make no mention of the special case $a_1=1$.  Some claims, like those in Example \ref{H&Vregions} below, should thus be understood up to this minor technicality, but this shall not affect the statements of any results.} of) the rectangles 
\[V_{a_2}\cap H_1,V_{a_2-1}\cap H_2,\dots,V_1\cap H_{a_2}\] 
and correspond to the Farey convergents 
\[\frac{p_{0}}{q_{0}},\frac{p_1+p_0}{q_1+q_0},\dots,\frac{(a_2-1)p_1+p_0}{(a_2-1)q_1+q_0},\]
respectively.  More generally, observe that $(x_n,y_n)\in H_1$ if and only if 
\[(x_n,y_n)=([0;a_{j+1},a_{j+2},\dots],[0;1,a_j,\dots,a_1])\] 
for some $j$ with $j>1$, or $j=1$ and $a_1>1$.  In that case, $n=a_1+a_2+\dots+a_j$ (i.e., $j=j_n(x)$; see Equation (\ref{j_lambda})), and the points
\[(x_n,y_n),(x_{n+1},y_{n+1}),\dots,(x_{n+a_{j+1}-1},y_{n+a_{j+1}-1})\]
belong to the rectangles
\[V_{a_{j+1}}\cap H_1,V_{a_{j+1}-1}\cap H_2,\dots,V_1\cap H_{a_{j+1}}\]
and correspond to the Farey convergents
\[\frac{p_{j-1}}{q_{j-1}},\frac{p_j+p_{j-1}}{q_j+q_{j-1}},\dots,\frac{(a_{j+1}-1)p_j+p_{j-1}}{(a_{j+1}-1)q_j+q_{j-1}},\]
respectively. 

\subsection{Inducing Ito's natural extension}\label{Inducing Ito's natural extension}

With the above identification of orbits and Farey convergents, we find that certain subregions $R\subset\Omega$ correspond to particular convergents or mediants of the {\sc rcf}-expansion of $x$; for instance, $H_{\lambda+1}$ corresponds to convergents ($\lambda=0$) or mediants ($\lambda>0$) of the form $(\lambda p_j+p_{j-1})/(\lambda q_j+q_{j-1}),\ \lambda<a_{j+1}$.  In particular, $H_1$ corresponds to {\sc rcf}-convergents $(p_j/q_j)_{j\ge -1}$, which are generated by the Gauss map $G$.  Moreover, the first return of $\F$ to $H_1$ is reminiscent of the natural extension $(\Omega,\mathcal{B},\bar\nu_G,\mathcal{G})$ of the Gauss map (cf. (\ref{GaussNE})):  if $(x,y)=([0;a_1,a_2,\dots],[0;b_1,b_2,\dots])\in H_1$, i.e., $b_1=1$, then 
\begin{align*}
\F^{a_1}([0;a_1,a_2,\dots],[0;1,b_2,b_3,\dots])=&\ \F^{a_1-1}([0;a_1-1,a_2,\dots],[0;2,b_2,b_3,\dots])\\
=&\ \F^{a_1-2}([0;a_1-2,a_2,\dots],[0;3,b_2,b_3,\dots])\\
&\vdots\\
=&\ \F([0;1,a_2,\dots],[0;a_1,b_2,b_3,\dots])\\
=&\ ([0;a_2,a_3\dots],[0;1,a_1,b_2,b_3,\dots]).
\end{align*}
In fact, Brown and Yin proved in \cite{BY1996} that a copy of the Gauss natural extension is found sitting (inverted, scaled and `suspended' from $y=1$) within $(\Omega,\mathcal{B},\bar\mu,\F)$:
\begin{thm}[Theorem 1 of \cite{BY1996}]\label{H_1&Gauss}
The Gauss natural extension $(\Omega,\mathcal{B},\bar\nu_G,\mathcal{G})$ is isomorphic via the map $(x,y)\mapsto (x,1/(y+1))$ to the dynamical system induced from $(\Omega,\mathcal{B},\bar\mu,\F)$ on the horizontal region $H_1$.
\end{thm}

These observations naturally lead one to consider $\F$ induced on other subregions $R\subset\Omega$ which pick out desired subsequences of Farey convergents.  

\begin{defn}
A $\bar\mu$-measurable subset $R\subset\Omega$ is called an \emph{inducible subregion} of $\Omega$ if either
\begin{enumerate}
\item[(i)] $R=\Omega$, or

\item[(ii)] $0<\bar\mu(R)<\infty$ and $R$ is a \emph{$\bar\mu$-continuity set}, i.e., $\bar\mu(\partial R)=0$.
\end{enumerate}
An inducible subregion $R$ satisfying (ii) is called \emph{proper}.
\end{defn}

\begin{remark}\label{ae_x_rem}
Due to Proposition \ref{Theta&h} below, our main interest is in $\F$-orbits of points of the form $(x,1)\in\Omega$ as in (\ref{(xn,yn)}).  The conditions of a proper inducible subregion $R$ guarantee that $\bar\mu(\text{int}(R))>0$ and hence for Lebesgue-a.e. $x\in[0,1]$, the $\F$-orbit of $(x,1)$ enters $R$ infinitely often.  (In fact, the stronger requirement that $\bar\mu(\partial R)=0$ is also needed for our purposes; see Remark \ref{mu_cont_needed}.)  Indeed, let $(x,z)\in \Omega$ with $x=[0;a_1,a_2,\dots]\notin\mathbb{Q}$.  For $n\ge a_1$, Equations (\ref{n_exp_with_x_digits}) and (\ref{F_NE_symbolic}) give that
\begin{equation}\label{F^n_explicit}
(x_n,z_n):=\F^n(x,z)=\begin{cases}
([0;a_{j_n+1}-\lambda_n,a_{j_n+2},\dots],[0;\lambda_n+1,a_{j_n},\dots,a_1-1+z^{-1}]),\ & z\neq 0,\\
([0;a_{j_n+1}-\lambda_n,a_{j_n+2},\dots],[0;\lambda_n+1,a_{j_n},\dots,a_2]),\ & z=0.
\end{cases}
\end{equation}
In particular, $z_n$ belongs to the cylinder of points in $[0,1]$ whose {\sc rcf}-expansions begin with $[0;\lambda_n+1,a_{j_n},\dots,a_2,\dots]$.  The Euclidean diameter of this cylinder is no greater than the reciprocal of the $j_{n}^\text{th}$ Fibonacci number squared and thus goes to $0$ uniformly in $z$ as $n$ goes to infinity (see also \cite{BY1996}).  Now let $E$ be the set of irrationals $x\in(0,1)$ for which $(x_n,y_n):=\F^n(x,1)$ enters $R$ at most finitely often.  Since $\bar\mu(\text{int}(R))>0$, there exists some point $(s,t)\in \text{int}(R)$ and some $\delta>0$ such that $B_{\delta}(s,t)\subset\text{int}(R)$.  If $E$ has positive Lebesgue measure, then also $\bar\mu(E\times[0,1])>0$.  Since $\F$ is conservative (\cite{BY1996}), $\bar\mu$-a.e.~point $(x,z)\in E\times [0,1]$ enters $B_{\delta/2}(s,t)$ infinitely often.  The observation on the diameter of the cylinder above implies that for $n$ large enough, $|(x_n,z_n)-(x_n,y_n)|<\delta/2$, and hence $(x_n,y_n)$ enters $B_{\delta}(s,t)\subset\text{int}(R)$ infinitely often---a contradiction. 
\end{remark}

For an inducible subregion $R\subset\Omega$, let $r=r_R:\Omega\to\mathbb{N}\cup\{\infty\}$ denote the hitting time
\begin{equation}\label{hitting_time}
r_R(x,y):=\inf\{\ n\ge 1\ |\ \F^n(x,y)\in R\}
\end{equation}
to $R$.  (Abusing notation, we assume that the null set of points in any set $S\subset \Omega$ which enter $R$ at most finitely many times under $\F$ is removed from $S$ and denote this new set again by $S$.)  Let $\F_R:\Omega\to R$ be defined by
\begin{equation}\label{ind_map}
\F_R(x,y):=\F^r(x,y)=\left(A_{[0,r]}^{-1}\cdot x\ ,\ A_{[r,0]}\cdot y\right), \quad \text{where} \quad r=r_R(x,y).
\end{equation}
The map $\F_R$ restricted to $R$ is the \emph{induced map} of $\F$ on $R$.  If $R$ is a proper inducible subregion, let $\bar\mu_R$ denote the \emph{induced measure}
\[\bar\mu_R(S):=\frac{\bar\mu(S)}{\bar\mu(R)},\quad S\in\mathcal{B}\cap R:=\{B\cap R\ |\ B\in\mathcal{B}\}.\] 
Since $(\Omega,\mathcal{B},\bar\mu,\F)$ is ergodic, so is the induced system $(R,\mathcal{B}\cap R,\bar\mu_R,\F_R)$.  In case $R=\Omega$ is not proper, we set $\bar\mu_R:=\bar\mu$ and note that $(\Omega,\mathcal{B},\bar\mu,\F)=(R,\mathcal{B}\cap R,\bar\mu_R,\F_R)$.  In this case, we may abuse terminology and refer to $\bar\mu_R$ as an induced measure and $(R,\mathcal{B}\cap R,\bar\mu_R,\F_R)$ as an induced system.  However, we emphasise that $\bar\mu_R$ is a finite (probability) measure if and only if $R$ is a \emph{proper} inducible subregion.

For an inducible subregion $R\subset \Omega$ and $(x,y)\in \Omega$, set $N_0^R(x,y):=0$ and, for $n\ge 1$,
\begin{equation}
N^R_n(x,y):=N_{n-1}^R(x,y)+r_R(\F_R^{n-1}(x,y)).
\end{equation}
In particular, this gives the following relationship between iterates of $\F_R$ and of $\F$, for each $n\ge 0$ and $(x,y)$:
\[\F_R^n(x,y)=\F^{N}(x,y)=\left(A_{[0,N]}^{-1}\cdot x\ ,\ A_{[N,0]}\cdot y\right), \quad \text{where} \quad N=N_n^R(x,y).\]
When the subregion $R$ and an irrational $x\in (0,1)$ are understood, we use the suppressed notation
\begin{equation}\label{N_n}
N_n:=N_n^R(x,1).
\end{equation}
An inducible subregion $R$ thus naturally determines a subsequence of the Farey convergents of Lebesgue-a.e. irrational\footnote{As suggested by Remark \ref{ae_x_rem}, there may be a Lebesgue-null set of points that enter $R$ at most finitely many times.  For instance, if $R=V_a$ as in Example \ref{H&Vregions} below and if $x=[0;a_1,a_2,\dots]$ with $a_j<a$ for all $j\ge 1$, then $\F^n(x,1)$ never enters $R$.  Such null sets are omitted from consideration throughout.} $x\in (0,1)$: for each $n\ge 0$, set $(x_n^R,y_n^R):=\F_R^n(x,1)$ and
\begin{equation}\label{induced_matrices}
\convmatR{n}=\begin{pmatrix}
u_n^R(x) & t_n^R(x)\\
s_n^R(x) & r_n^R(x)
\end{pmatrix}
:=A_{[0,N_n]}=\convmat{N_n}=\begin{pmatrix}\lambda p_j+p_{j-1} & p_j\\ \lambda q_j+q_{j-1} & q_j\end{pmatrix},
\end{equation}
where $j=j_{N_n}(x)$ and $\lambda=\lambda_{N_n}(x)$ (recall (\ref{j_lambda}) and (\ref{convmats})).
Informally speaking, the subsequence 
\[(u_n^R/s_n^R)_{n\ge 0}
=(u_{N_n}/s_{N_n})_{n\ge 0}=(P_{N_n-1}/Q_{N_n-1})_{n\ge0}\]
of Farey convergents corresponding to $R$ consists of those convergents which are `picked up' when the forward orbit of $(x,1)$ under $\F$ enters the region $R$.

\begin{eg}\label{H&Vregions}
As noted above, the region $R=H_1$ corresponds to the {\sc rcf}-convergents of $x$ (see Figure \ref{MetricFigs}.i).  In particular, 
\[(u_n^{H_1}/s_n^{H_1})_{n\ge 0}=(p_{j-1}/q_{j-1})_{j\ge 0}.\]
Moreover, for $\lambda\ge 1$ the region $R=H_{\lambda+1}$ gives the $\lambda^{\text{th}}$ mediant convergents
\[\{u_n^{H_{\lambda+1}}/s_n^{H_{\lambda+1}}\}_{n\ge 0}=\{(\lambda p_j+p_{j-1})/(\lambda q_j+q_{j-1})\ |\ \lambda<a_{j+1}\}_{j\ge 0}\]
(Figure \ref{MetricFigs}.ii).  Similarly, the vertical regions $V_a,\ a=1,2,\dots,$ give---in addition to the {\sc rcf}-convergents $p_{j-1}/q_{j-1}$ for which $a_{j+1}=a$---the final mediants, next-to-final mediants, and so on, respectively (Figure \ref{MetricFigs}.iii): 
\[\{u_n^{V_{a}}/s_n^{V_a}\}_{n\ge 0}=\{((a_{j+1}-a) p_j+p_{j-1})/((a_{j+1}-a) q_j+q_{j-1})\ |\ a_{j+1}\ge a\}_{j\ge 0}.\]
\end{eg}

\begin{eg}\label{rectangles_eg}
The rectangles $R=V_a\cap H_b$ pick out particular convergents ($b=1$) or mediant convergents ($b>1$) corresponding to specific partial quotients in the {\sc rcf}-expansion of $x$.  For instance, 
\[\{u_n^{V_3\cap H_1}/s_n^{V_3\cap H_1}\}_{n\ge 0}=\{p_{j-1}/q_{j-1}\ |\ a_{j+1}=3\}_{j\ge0}\]
\[\{u_n^{V_2\cap H_2}/s_n^{V_2\cap H_2}\}_{n\ge 0}=\{(p_j+p_{j-1})/(q_j+q_{j-1})\ |\ a_{j+1}=3\}_{j\ge 0},\]
and
\[\{u_n^{V_1\cap H_3}/s_n^{V_1\cap H_3}\}_{n\ge 0}=\{(2p_j+p_{j-1})/(2q_j+q_{j-1})\ |\ a_{j+1}=3\}_{j\ge 0}\]
(see the three left-most plots of Figure \ref{FareyNE}).

More generally, for $R=V_{a-\lambda}\cap H_{\lambda+1}$ with $a>0$ and $0\le \lambda<a$,
\[\{u_n^R/s_n^R\}_{n\ge 0}=\{(\lambda p_j+p_{j-1})/(\lambda q_j+q_{j-1})\ |\ a_{j+1}=a\}_{j\ge 0}.\]
\end{eg}

We end this subsection with a calculation of the measure-theoretic entropy $h(\F_R)=h_{\bar\mu_R}(\F_R)$ of the induced transformation $\F_R$ restricted to a proper inducible $R\subset\Omega$.

\begin{thm}\label{entropy_thm}
Let $R\subset\Omega$ be a proper inducible subregion.  Then 
\[h(\F_R)=\frac{\pi^2}{6\bar\mu(R)}.\]
\end{thm}
\begin{proof}
We first note that for two proper inducible subregions $R_1,R_2\subset \Omega$ with $R_1\subset R_2$, the dynamical system $(R_1,\mathcal{B}\cap R_1,\bar\mu_{R_1},\F_{R_1})$ is isomorphic to \emph{the induced system of $(R_2,\mathcal{B}\cap R_2,\bar\mu_{R_2},\F_{R_2})$ on $R_1$}.  Hence, by Abramov's formula,
\[
h(\F_{R_1})=\frac{h(\F_{R_2})}{\bar\mu_{R_2}(R_1)}=\frac{\bar\mu(R_2)}{\bar\mu(R_1)}h(\F_{R_2}),
\]
or
\begin{equation}\label{abramov}
\bar\mu(R_1)h(\F_{R_1})=\bar\mu(R_2)h(\F_{R_2}).
\end{equation}
It is well-known that the entropy of the Gauss map $G$ (and hence also its natural extension $\mathcal{G}$) is $\pi^2/6\log2$, and thus $h(\F_{H_1})=\frac{\pi^2}{6\log2}$ by Theorem \ref{H_1&Gauss}.  Using this, (\ref{abramov}) and a calculation of $\bar\mu(H_1)=\log2$, we compute
\[\bar\mu(R)h(\F_R)=\bar\mu(R\cup H_1)h(\F_{R\cup H_1})=\bar\mu(H_1)h(\F_{H_1})=\frac{\pi^2}{6}.\]
\end{proof}

\subsection{The (relative) equidistribution of ($\F$-) $\F_R$-orbits}\label{The (relative) equidistribution of orbits}

In \cite{BY1996}, Brown and Yin employ the Ratio Ergodic Theorem to derive metrical results on the system $(\Omega,\mathcal{B},\bar\mu,\F)$:
\begin{thm}[Theorem 2 of \cite{BY1996}]\label{RatThm}
For any $f,g\in L^1(\bar\mu)$ with $\int gd\bar\mu\neq 0$, 
\begin{equation}\label{RET}
\lim_{n\to\infty}\frac{\sum_{k=0}^{n-1}f(x_k,z_k)}{\sum_{k=0}^{n-1}g(x_k,z_k)}=\frac{\int fd\bar\mu}{\int gd\bar\mu} \quad \text{$\bar\mu$-a.s.,}
\end{equation}
where $(x_k,z_k):=\F^k(x,z),\ k\ge 0$.
\end{thm}
Under certain Lipschitz conditions on $f$ and $g$, Brown and Yin are able to replace $(x_k,z_k)$ on the left-hand side of (\ref{RET}) with $(x_k,y_k):=\F^k(x,1)$ for almost every\footnote{All \emph{almost every} statements are w.r.t.\ Lebesgue measure.}  $x\in(0,1)\backslash \mathbb{Q}$ (see Theorem 3 of \cite{BY1996}).  The following theorem---which does not require Lipschitz conditions but instead replaces $f$ and $g$ by indicator functions---may be seen as an analogue of an important result of Jager for the natural extension of the Gauss map, which states that for almost every irrational $x\in(0,1)$, the $\mathcal{G}$-orbit of $(x,0)$ is $\bar\nu_G$-equidistributed over $\Omega$ (see Theorem 3 of \cite{J86}).  In fact, by Theorem \ref{H_1&Gauss} above and Corollary \ref{muR_equidist} below, the following statement may be read as a generalisation of Jager's result:
\begin{thm}\label{rel-equidist}
For almost every $x\in(0,1)\backslash\mathbb{Q}$, the $\F$-orbit of $(x,1)$ is \emph{$\bar\mu$-relatively equidistributed}\footnote{This notion of $\bar\mu$-relative equidistribution differs slightly from that studied by Gerl in \cite{G1971A} and \cite{G1971B}.  In particular, Gerl required a Radon measure, but $\bar\mu$ is not locally finite.} over $\Omega$.  That is, for almost all $x\in(0,1)\backslash\mathbb{Q},$
\[\lim_{n\to\infty}\frac{\sum_{k=0}^{n-1}\mathbf{1}_S(x_k,y_k)}{\sum_{k=0}^{n-1}\mathbf{1}_R(x_k,y_k)}=\frac{\bar\mu(S)}{\bar\mu(R)}\]
for any proper inducible or $\bar\mu$-null set $S\subset\Omega$ and any proper inducible $R\subset\Omega$.
\end{thm}
\begin{proof}
Let $\mathcal{U}$ be a countable base for the subspace topology on $\Omega\subset \mathbb{R}^2$, and let $\mathcal{C}$ denote the countable collection of all finite unions of finite $\bar\mu$-measure elements of $\mathcal{U}$, along with the empty set.  Moreover, for $A,B\in\mathcal{C}$ with $\bar\mu(B)>0$, let $E(A,B)\subset (0,1)\backslash\mathbb{Q}$ denote the set of irrational $x$ for which there is \emph{no} $z\in [0,1]$ satisfying 
\begin{equation}\label{RETC}
\lim_{n\to\infty}\frac{\sum_{k=0}^{n-1}{\mathbf{1}_{A}(x_k,z_k)}}{\sum_{k=0}^{n-1}{\mathbf{1}_{B}(x_k,z_k)}}=\frac{\bar\mu(A)}{\bar\mu(B)}.
\end{equation}
Theorem \ref{RatThm} implies that 
$E(A,B)$ is a Lebesgue-null set: otherwise $E(A,B)\times[0,1]$ is a set of positive $\bar\mu$-measure for which (\ref{RET}) does not hold.  Since $\mathcal{C}$ is countable, the union $E$ of all such $E(A,B)$ is also a Lebesgue-null set.  Hence, for every $x\in (0,1)\backslash (E\cup \mathbb{Q})$ it follows that for all $A, B\in\mathcal{C}$ with $\bar\mu(B)>0$, there is some $z\in[0,1]$ for which (\ref{RETC}) holds.

We claim that for any proper inducible or $\bar\mu$-null set $S\subset\Omega$ and any $\delta>0$, there exist $S_{+\delta},S_{-\delta}\in\mathcal{C}$ such that 
\begin{enumerate}
\item[(i)] $S_{-\delta}\subset S\subset S_{+\delta},$

\item[(ii)] $\bar\mu(S_{+\delta}\backslash S),\ \bar\mu(S\backslash S_{-\delta})<\delta$, and

\item[(iii)] $d(S, \Omega\backslash S_{+\delta})>0$, and if $\bar\mu(S)>0$, also $d(S_{-\delta}, \Omega\backslash S)>0$, 
\end{enumerate}
where for $A,B\in\mathcal{B}$, 
\[d(A,B):=\inf\left\{|a-b|\ \big|\ a\in A,\ b\in B\right\}\]
denotes the Euclidean distance between the sets $A$ and $B$.  Indeed, fix $S\subset\Omega$ and $\delta>0$ as above.  By the regularity of $\bar\mu$, there exists an open cover $\{U_i\}_{i\in I}\subset \mathcal{U}$ of the closure $\bar S$ of $S$ for which $\bar\mu(\cup_{i\in I}U_i\backslash \bar S)<\delta$.  Since $\bar S$ is compact, $\{U_i\}_{i\in I}$ has some finite subcover, the union of whose elements we denote by $S_{+\delta}$.  Note, then, that $S\subset S_{+\delta}\subset \cup_{i\in I}U_i$, and since the boundary $\partial S$ is a $\bar\mu$-null set, $\bar\mu(S_{+\delta}\backslash S)\le \bar\mu(\cup_{i\in I}U_i\backslash \bar S)<\delta$.  Moreover, since $\bar S$ and $\Omega\backslash S_{+\delta}$ are compact and disjoint, the distance between them is strictly positive, and thus the distance between $S\subset \bar S$ and $\Omega\backslash S_{+\delta}$ is strictly positive.  Thus $S_{+\delta}$ satisfies each of the properties of the claim. 

If $\bar\mu(S)=0$, set $S_{-\delta}:=\varnothing\in\mathcal{C}$, which trivially satisfies the claim.  Now suppose $\bar\mu(S)>0$.  Since $\bar\mu(\partial S)=0$, we have $\bar\mu(\text{int}(S))>0$.  Again by regularity of $\bar\mu$, there exists some compact subset $K$ of $\text{int}(S)$ with $\bar\mu(\text{int}(S)\backslash K)<\delta$.  Since $\Omega$ is normal and $\Omega\backslash \text{int}(S)$ and $K$ are closed and disjoint, there exist open, disjoint sets $U,\ V\in\mathcal{B}$ containing $\Omega\backslash \text{int}(S)$ and $K$, respectively.  Let $\{V_j\}_{j\in J}$ be a collection of open sets from the countable base $\mathcal{U}$ whose union equals $V$.  This collection of sets forms an open cover of the compact set $K$; let $S_{-\delta}\in\mathcal{C}$ be the union of elements of a finite subcover.  We then have $S_{-\delta}\subset \text{int}(S)\subset S$ and, since $\bar\mu(\partial S)=0$, also $\bar\mu(S\backslash S_{-\delta})\le\bar\mu(\text{int}(S)\backslash K)<\delta$.  Furthermore, since $S_{-\delta}\subset \Omega\backslash U$ and $\Omega\backslash S\subset \Omega\backslash \text{int}(S)$, we have $d(S_{-\delta},\Omega\backslash S)\ge d(\Omega\backslash U,\Omega\backslash \text{int}(S))$.  But $\Omega\backslash U$ and $\Omega\backslash \text{int}(S)$ are compact and disjoint, so the distance between them is again strictly positive.  Thus $S_{-\delta}$ also satisfies the desired properties. 

For $S\subset\Omega$ and $\delta>0$ as above, let
\[
d_S(\delta):=\begin{cases}
d(S, \Omega\backslash S_{+\delta}),\ & \bar\mu(S)=0,\\
\min\{d(S, \Omega\backslash S_{+\delta}),d(S_{-\delta}, \Omega\backslash S)\},\ & \bar\mu(S)>0.
\end{cases}
\]
Recall from Remark \ref{ae_x_rem} that for any $x\in(0,1)\backslash\mathbb{Q}$, the distances between $(x_n,y_n):=\F^n(x,1)$ and $(x_n,z_n):=\F^n(x,z)$ approach zero uniformly in $z\in[0,1]$.  Thus for any $x\in(0,1)\backslash\mathbb{Q}$ there exists some $n_S(\delta)\in\mathbb{N}$ such that for any $z\in[0,1]$ and $n\ge n_S(\delta)$,
\[|(x_n,y_n)-(x_n,z_n)|<d_S(\delta).\]
In particular, by definition of $d_S(\delta)$, for any $n\ge n_S(\delta)$, 
\[(x_n,y_n)\in S \quad \text{implies} \quad (x_n,z_n)\in S_{+\delta},\]
and if $\bar\mu(S)>0$,
\[(x_n,z_n)\in S_{-\delta} \quad \text{implies} \quad (x_n,y_n)\in S.\]

Now let $x\in(0,1)\backslash(E\cup\mathbb{Q})$, $S\subset\Omega$ a proper inducible subregion or $\bar\mu$-null set and $R\subset\Omega$ a proper inducible subregion.  For $\delta>0$, let $S_{\pm\delta},R_{\pm\delta}\in\mathcal{C}$ be the sets constructed above.  When $\delta$ is sufficiently small, the observations at the beginning of the proof imply that there exist $z,z'\in[0,1]$ such that
\begin{align*}
\frac{\bar\mu(S_{-\delta})}{\bar\mu(R_{+\delta})}=\liminf_{n\to\infty}\frac{\sum_{k=0}^{n-1}{\mathbf{1}_{S_{-\delta}}(x_k,z_k)}}{\sum_{k=0}^{n-1}{\mathbf{1}_{R_{+\delta}}(x_k,z_k)}}\le&\liminf_{n\to\infty}\frac{\sum_{k=0}^{n-1}{\mathbf{1}_{S}(x_k,y_k)}}{\sum_{k=0}^{n-1}{\mathbf{1}_{R}(x_k,y_k)}}\\
\le &\limsup_{n\to\infty}\frac{\sum_{k=0}^{n-1}{\mathbf{1}_{S}(x_k,y_k)}}{\sum_{k=0}^{n-1}{\mathbf{1}_{R}(x_k,y_k)}}\\
\le &\limsup_{n\to\infty}\frac{\sum_{k=0}^{n-1}{\mathbf{1}_{S_{+\delta}}(x_k,z'_k)}}{\sum_{k=0}^{n-1}{\mathbf{1}_{R_{-\delta}}(x_k,z'_k)}}=\frac{\bar\mu(S_{+\delta})}{\bar\mu(R_{-\delta})}.
\end{align*}
By properties (i) and (ii) above, taking $\delta\to 0$ gives the result.  
\end{proof}

As a corollary, we obtain the following:
\begin{cor}\label{muR_equidist}
For almost every $x\in (0,1)\backslash\mathbb{Q}$, the $\F_R$-orbit of $(x,1)$ is $\bar\mu_R$-equidistributed for any proper inducible $R\subset\Omega$.  That is, for almost every $x\in (0,1)\backslash\mathbb{Q}$,
\[\lim_{n\to\infty}\frac1n\sum_{k=0}^{n-1}\mathbf{1}_S(x_k^R,y_k^R)=\bar\mu_R(S)\]
for any proper inducible $R\subset\Omega$ and any $S\in\mathcal{B}\cap R$ with $\bar\mu_R(\partial S)=0$.
\end{cor}
\begin{proof}
When $S\subset R$, 
\[\lim_{n\to\infty}\frac1n\sum_{k=0}^{n-1}\mathbf{1}_S(x_k^R,y_k^R)=\lim_{n\to\infty}\frac{\sum_{k=0}^{n-1}\mathbf{1}_S(x_k^R,y_k^R)}{\sum_{k=0}^{n-1}\mathbf{1}_R(x_k^R,y_k^R)}=\lim_{n\to\infty}\frac{\sum_{k=0}^{N_n-1}\mathbf{1}_S(x_k,y_k)}{\sum_{k=0}^{N_n-1}\mathbf{1}_R(x_k,y_k)},\]
where $N_n=N_n^R(x,1)$.  The result follows from Theorem \ref{rel-equidist} and the fact that $\bar\mu_R(S)=\bar\mu(S)/\bar\mu(R)$.
\end{proof}

\begin{remark}\label{mu_cont_needed}
The $\bar\mu$-continuity condition on $R$, namely that $\bar\mu(\partial R)=0$, is necessary to avoid `pathological' counter-examples to the above result.  Indeed, since $y_n^R\in\mathbb{Q}$ for all $n$, if $R\subset[0,1]\times ([0,1]\backslash\mathbb{Q})$, then $(x_n^R,y_n^R)$ never enters $R$.  Using similar ideas, one can easily construct $R\in\mathcal B$ with $0<\bar\mu(R)<\infty$ and $\bar\mu(\partial R),\ \bar\mu(\text{int}(R))>0$ so that $(x_n^R,y_n^R)$ almost surely enters $R$ infinitely often but the conclusion of Corollary \ref{muR_equidist} is false (cf. Remark \ref{ae_x_rem}).
\end{remark}

\bigskip
\section{Metrical results}\label{Metrical results}

An important reason to study the {\sc rcf}-expansion is that this algorithm yields rational approximations to irrational numbers of `very high quality.' By this, we mean that the {\sc rcf}-convergents $(p_n/q_n)_{n\geq 0}$ of an irrational number $x$ have the so-called \emph{best approximation property} (recall (\ref{best_approx_ineq})). This result was essentially already known to Christiaan Huygens when he constructed his planetarium in 1680 at the request of Jean-Baptiste Colbert (see Chapter IV in \cite{RS92}). Later in this section (\S\ref{On the theorems of Legendre, Fatou--Grace and Koksma}) we will revisit an old result by Legendre and some of its refinements. Legendre's result states that if $p/q$ is a rational number in its lowest terms with $q$ positive, and
\begin{equation}\label{Legendre1}
\left| x - \frac{p}{q}\right| < \frac{1}{2q^2},
\end{equation}
one has that there exists an $n\in\mathbb{N}$ such that $p=p_n$ and $q=q_n$. In other words, in order to approximate an irrational $x$ `well' by a rational $p/q$ (`well' in the sense that (\ref{Legendre1}) holds) one is bound to find a {\sc rcf}-convergent of $x$. In view of this, for over a century so-called \emph{{\sc rcf}-approximation coefficients} $\theta_n$, defined by
\[\theta_n=\theta_n(x) := q_n^2\left| x - \frac{p_n}{q_n}\right| ,\qquad  \text{for $n\in\mathbb{N}$},\]
have been studied.  Independently, Doeblin and Lenstra conjectured the distribution for almost all $x$ of the sequence $(\theta_n(x))_{n\in\mathbb{N}}$. The proof of this so-called Doeblin--Lenstra conjecture by Bosma, Jager and Wiedijk in \cite{BJW83} (cf. the approach of Knuth in \cite{K1984}) lead to many new results in Diophantine approximation, some of which are mentioned below.  In \S\ref{Approximation coefficients and their limiting distributions} we consider a number of old and new Doeblin--Lenstra-type theorems for subsequences of approximation coefficients corresponding to {\sc rcf}-convergents and mediants.

Apart from the Legendre-type results considered in \S\ref{On the theorems of Legendre, Fatou--Grace and Koksma}, we mention here some classical results which hold for all irrational $x$ and all $n\in\mathbb{N}$:
\[\min \{ \theta_{n-1}(x),\theta_n(x)\} < \frac{1}{2} \quad \text{(Vahlen, 1895 \cite{V95}; see also \S\ref{Consecutive rcf-convergents} below)};\]
\[\min \{ \theta_{n-1}(x),\theta_n(x),\theta_{n+1}(x)\} \leq \frac{1}{\sqrt{5}} \quad \text{(Borel, 1903 \cite{B03})};\]
which in itself is a corollary of the following result
\[\min \{ \theta_{n-1}(x),\theta_n(x),\theta_{n+1}(x)\} < \frac{1}{\sqrt{a_{n+1}^2+4}} \quad \text{(Bagemihl \& McLaughlin, 1966 \cite{BM66})}.\]
Related to this last result is the result by Tong from 1983 (\cite{T83}), which states that:
\[\max \{ \theta_{n-1}(x),\theta_n(x),\theta_{n+1}(x)\} > \frac{1}{\sqrt{a_{n+1}^2+4}}.\]
These results can easily be derived using the natural extension of the Gauss map (see e.g.~\cite{JK89,DK2002B}), and thus by Theorem \ref{H_1&Gauss} may also be derived through the set-up of this paper.  Subsection \ref{Consecutive approximation coefficients} builds a general framework for studying consecutive approximation coefficients corresponding to subsequences of {\sc rcf}-convergents and mediants determined by inducible $R\subset \Omega$.

In 1936, Paul L{\'e}vy (\cite{L36}) proved the following important and classical result: for almost all $x\in[0,1]$,
\begin{equation}\label{levy}
\lim_{n\to\infty}\frac1n\log q_n=\frac{\pi^2}{12\log 2}.
\end{equation}
As a corollary of this and the fact that $1/2q_nq_{n+1}<|x-p_n/q_n|<1/q_nq_{n+1}$, it follows that for almost every $x\in[0,1]$,
\begin{equation}\label{levy2}
\lim_{n\to\infty}\frac1n\log\left|x-\frac{p_n}{q_n}\right|=-\frac{\pi^2}{6\log2}.
\end{equation}
In \S\ref{A generalised Levy-type theorem}, we obtain a L\'evy-type theorem for subsequences of {\sc rcf}-convergents and mediants corresponding to proper inducible subregions $R\subset \Omega$, which generalises a number of existing results from the literature.  It should be mentioned that originally L\'evy's result (and similar other results by L\'evy and Khintchine) also gave a `speed of convergence,' as its proof relied on probability theory and was derived from the Gauss-Kuzmin-L\'evy Theorem; see e.g.~\cite{IK02} for more details. As we use ergodic theory, such speeds cannot be given.

\subsection{Approximation coefficients and their limiting distributions}\label{Approximation coefficients and their limiting distributions}
In this subsection we obtain---among new results---several metrical theorems of \cite{BJW83,I1989,B1990,J91,BY1996} as simple corollaries.  Some of the theorems of \cite{B1990} are obtained in a similar fashion in \cite{IK2008} using the natural extension of \emph{Denjoy's canonical continued fraction map} $T_d:[0,\infty)\to [0,\infty)$, defined by
\[
T_d(x):=\begin{cases}
0 & x=0,\\
\frac{1-x}{x} & x\le 1,\\
\frac1{x} & x>1.
\end{cases}
\]
In fact, the Farey tent map $F$ is the first-return map of $T_d$ to $[0,1]$, and the proofs of \cite{IK2008} and those found below are closely related.  However, while the domain $\Omega$ of the natural extension of $F$ is bounded, the domain of the natural extension of $T_d$ considered in \cite{IK2008} is the unbounded region $([0,1)\times[0,\infty))\cup([1,\infty)\times[0,1])$.  We find the new proofs in the setting of $(\Omega,\mathcal{B},\bar\mu,\F)$ to be particularly insightful given their concrete geometric realisation within $\Omega$.

For any $x\in\mathbb{R}$ and $p/q\in\mathbb{Q}$ with $\text{gcd}\{p,q\}=1$ and $q>0$, set
\[\Theta(x,p/q):=q|qx-p|.\]
Observe that 
\[\left|x-\frac{p}{q}\right|=\frac{\Theta(x,p/q)}{q^2},\]
so the \textit{approximation coefficient} $\Theta(x,p/q)$ gives a measure of how well the rational $p/q$ approximates $x$.  Notice that $\Theta(x,p_n/q_n)=\theta_n(x)$ when $p_n/q_n$ is the $n^\text{th}$ {\sc rcf}-convergent of $x$.  When $p/q=u_n/s_n=P_{n-1}/Q_{n-1}$ is the $(n-1)^\text{st}$ Farey convergent of $x$, we use the special notation
\[\Theta_n(x):=\Theta(x,u_n/s_n).\]
Similarly, for an inducible subregion $R\subset\Omega$, we set
\begin{equation}\label{Theta_n^R}
\Theta_n^R(x):=\Theta(x,u_n^R/s_n^R).
\end{equation}

The following result---which is central to the remainder of the paper---states that the $n^\text{th}$ approximation coefficient $\Theta_n(x)$ is computed explicitly in terms of the $n^\text{th}$ point $(x_n,y_n)$ in the $\F$-orbit of $(x,1)$ and thus depends on both the `future' and the `past' of the $F$-orbit of $x$.  Hence---although the statements of results regarding approximations coefficients are about $x\in (0,1)$---the proofs exploit the two-dimensional system $(\Omega,\mathcal{B},\bar\mu,\F)$.  Define $h:\Omega\backslash\{(0,0)\}\to[0,\infty)$ by 
\[h(x,y):=\frac{1-y}{x+y-xy}.\]

\begin{prop}[cf. Propositions 1.2 and 2.2 of \cite{I1989}]\label{Theta&h}
For any inducible $R\subset\Omega$ and $n\ge 0$,
\[\Theta_n(x)=h(x_n,y_n)\quad \text{and}\quad \Theta_n^R(x)=h(x_n^R,y_n^R).\]
\end{prop}
\begin{proof}
Notice from Equations (\ref{convmats_reversed}) and (\ref{(xn,yn)}) that 
\[y_n=A_{[n,0]}\cdot 1=\frac{r_n}{s_n+r_n},\]
so
\[h(x_n,y_n)=\frac{1-y_n}{x_n+y_n-x_ny_n}=\frac{1-y_n}{(1-y_n)x_n+y_n}=\frac{s_n}{s_nx_n+r_n}.\]
On the other hand, 
\[x=A_{[0,n]}\cdot x_n=\frac{u_n x_n+t_n}{s_n x_n+r_n},\]
so
\[\Theta_n(x)=s_n^2\left|x-\frac{u_n}{s_n}\right|=s_n^2\left|\frac{u_n x_n+t_n}{s_n x_n+r_n}-\frac{u_n}{s_n}\right|=\frac{s_n\left|t_ns_n-u_nr_n\right|}{s_n x_n+r_n}=\frac{s_n}{s_n x_n+r_n}=h(x_n,y_n),\]
where the penultimate equality follows from $\text{det}(A_{[0,n]})=\pm1$.  The second equality of the proposition statement now follows from the first:
\[\Theta_n^R(x)=\Theta(x,u_n^R/s_n^R)=\Theta(x,u_{N_n}/s_{N_n})=\Theta_{N_n}(x)=h(x_{N_n},y_{N_n})=h(x_n^R,y_n^R).\]
\end{proof}

\begin{figure}[t]
\includestandalone[width=.4\textwidth]{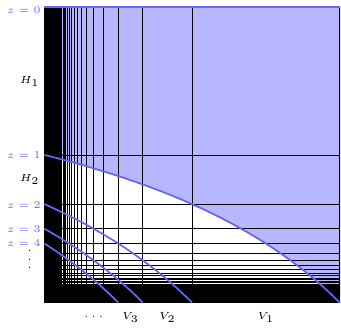}
\caption{The curves $h(x,y)=z$ for $z\in\{0,1,2,3,4\}$.  The region $S_1$ is the shaded region above the curve $h(x,y)=1$.}
\label{SzFig}
\end{figure}

For $z\in[0,\infty)$, let
\[S_z:=\{(x,y)\in\Omega\ |\ h(x,y)\le z\};\]
see Figure \ref{SzFig}.  The previous result together with Corollary \ref{muR_equidist} allows us to calculate the asymptotic relative frequency of bounded approximation coefficients $\Theta_n^R(x)\le z$ as the $\bar\mu_R$-measure of $S_z\cap R$:

\begin{thm}\label{approx_coeff_freq}
For any proper inducible subregion $R\subset\Omega$, almost every $x\in (0,1)\backslash\mathbb{Q}$ satisfies 
\[\lim_{n\to\infty}\frac1n\#\{0\le k< n\ |\ \Theta_k^R(x)\le z\}=\bar\mu_R(S_z\cap R) \quad \text{for all} \quad z\in[0,\infty).\]
\end{thm}
\begin{proof}
With the stated assumptions, Proposition \ref{Theta&h} and Corollary \ref{muR_equidist} give that for almost every irrational $x\in(0,1)$, it follows for every $z\in[0,\infty)$ 
\begin{align*}
\lim_{n\to\infty}\frac1n\#\{0\le k< n\ |\ \Theta_k^R(x)\le z\}=&\lim_{n\to\infty}\frac1n\#\{0\le k< n\ |\ h(x_k^R,y_k^R)\le z\}\\
=&\lim_{n\to\infty}\frac1n \sum_{k=0}^{n-1}\mathbf{1}_{S_z\cap R}(x_k^R,y_k^R)\\
=&\bar\mu_R(S_z\cap R).
\end{align*}
\end{proof}

\begin{figure}[t]
\includestandalone[width=.22\textwidth]{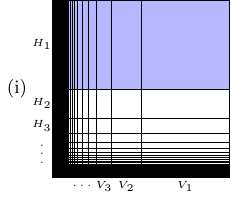}
\includestandalone[width=.22\textwidth]{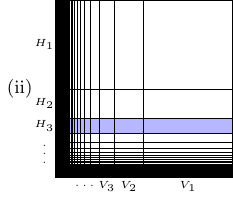}
\includestandalone[width=.22\textwidth]{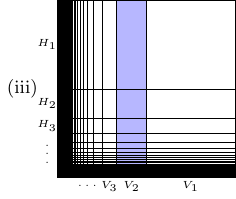}
\includestandalone[width=.22\textwidth]{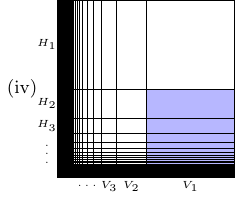}
\includestandalone[width=.22\textwidth]{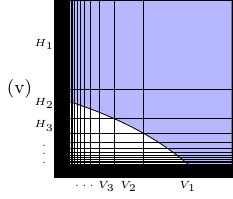}
\includestandalone[width=.22\textwidth]{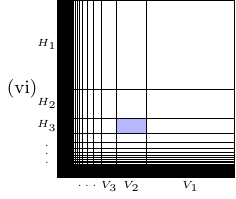}
\includestandalone[width=.22\textwidth]{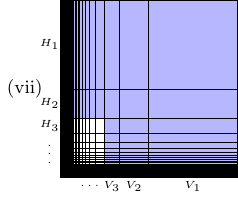}
\caption{The regions $R$ considered in Corollaries \ref{cor.i}--\ref{cor.vi} and \ref{cor.vii}.}
\label{MetricFigs}
\end{figure}

Theorem \ref{approx_coeff_freq} allows us to easily obtain a number of classical and new results on the limiting distributions of approximation coefficients corresponding to particular subsequences of {\sc rcf}-convergents and mediants of generic $x\in(0,1)\backslash\mathbb{Q}$.  The proofs of the corollaries below follow the same general procedure.  For a given proper inducible subregion $R\subset\Omega$, one first finds the infimum, $z_-$, and supremum, $z_+$, values of $z$ for which $h(x,y)=z$ intersects $R$.  By Proposition \ref{Theta&h}, these give bounds on $\Theta_n^R(x)$.  That the bounds are optimal follows from Theorem \ref{approx_coeff_freq}: for any $z\in (z_-,z_+)$ one finds $\bar\mu_R(S_z\cap R)\in (0,1)$, and thus the limiting distribution of Theorem \ref{approx_coeff_freq} implies that for almost every $x\in (0,1)\backslash\mathbb{Q}$, there exist $n$ for which $\Theta_n^R(x)\le z$ and $n$ for which $\Theta_n^R(x)>z$.  Moreover, Theorem \ref{approx_coeff_freq} leads to an explicit computation of the asymptotic relative frequency
\[\lim_{n\to\infty}\frac1n\#\{0\le k< n\ |\ \Theta_k^R(x)\le z\};\]
it is simply a matter of calculating the $\bar\mu$-measures of $S_z\cap R$ and $R$ and taking their quotient.  The former measure depends on how the curve $h(x,y)=z$ intersects the region $R$, and these different possible intersections give rise to the different branches of the frequency functions seen below.  The limiting distributions in each of the following corollaries are plotted here (\url{https://www.desmos.com/calculator/pnsu8hcytq}) via Desmos (\cite{desmos}).  

\begin{remark}\label{rearr_remark}
Each of the limiting distributions that we `re-obtain' below were originally formulated assuming that the corresponding subsequences of {\sc rcf}-convergents and mediants are ordered with increasing denominators.  However, recall from (\ref{FLseq}) that the denominators of the sequence $(u_n/s_n)_{n\ge 0}$ of Farey convergents are not necessarily increasing, and hence the same is true of the subsequence $(u_n^R/s_n^R)_{n\ge 0}$ for inducible $R\subset \Omega$.  Proposition \ref{rearr_1}---the statement and proof of which are found in the appendix (\S\ref{Appendix})---states that for the purposes of these limiting distributions, our reordering is innocuous and thus the original results do in fact follow from our methods.
\end{remark}

The first of our corollaries is the Doeblin--Lenstra conjecture, which was first proven in \cite{BJW83}; see Example \ref{H&Vregions} and Figure \ref{MetricFigs}.i.  For this we sketch the computations outlined above Remark \ref{rearr_remark}; proofs of the remaining corollaries are similar, though at times more tedious.

\begin{cor}[Doeblin--Lenstra, Theorem 1 of \cite{BJW83}]\label{cor.i}
Set $R:=H_1$ so that for any $x=[0;a_1,a_2,\dots]\in(0,1)\backslash \mathbb{Q}$,
\[\{u_n^R/s_n^R\}_{n\ge 0}=\{p_{j-1}/q_{j-1}\}_{j\ge 0}.\]
Then for $n>0$,
\[0<\Theta_n^R(x)<1,\]
with the upper and lower bounds optimal, and for almost every such $x$,
\[
\lim_{n\to\infty}\frac1n\#\{0\le k< n\ |\ \Theta_k^R(x)\le z\}
=\begin{cases}
\frac1Cz, & 0\le z\le 1/2,\\
\frac1C\left(1-z+\log\left(2z\right)\right), & 1/2\le z\le 1,
\end{cases}
\]
where $C:=\bar\mu(H_1)=\log 2$.
\end{cor}
\begin{proof}
The lower and upper bounds on $\Theta_n^R(x)$ follow from the argument above Remark \ref{rearr_remark}: the infimum, $z_-$, and supremum, $z_+$, of values $z$ for which $h(x,y)=z$ intersects $R=H_1$ are $0$ and $1$, respectively (see Figure \ref{SzFig}).  

For the limiting distribution of approximation coefficients, it suffices to compute $\bar\mu_R(S_z\cap H_1)$ for each $z\in[0,1]$, and to set $C=\bar\mu(H_1)=\bar\mu(S_1\cap H_1)$.  Fix $z\in[0,1]$.  Rearranging terms, $h(x,y)=z$ may be written as $y=f(x,z)$, where
\[f(x,z)=\frac{1-xz}{1-xz+z}.\]
Now $\frac{\partial f}{\partial x}\le 0,$ so $f$ is monotone decreasing as a function of $x$.  For $0\le z\le 1/2$, both $f(0,z)=1/(1+z)$ and $f(1,z)=1-z$ are bounded between $1/2$ and $1$, so the same is true of $f(x,z)$ for all $0\le x\le 1$.  Hence for $0\le z\le 1/2,$ one calculates
\[\bar\mu(S_z\cap H_1)=\int_0^1\int_{f(x,z)}^1\frac{dydx}{(x+y-xy)^2}=z.\]
On the other hand, if $1/2\le z\le 1$, then $1/2\le f(0,z)\le 1$ while $f(1,z)\le 1/2$.  One finds that $f((1-z)/z,z)=1/2$ and thus computes
\[\bar\mu(S_z\cap H_1)=\int_0^{\frac{1-z}{z}}\int_{f(x,z)}^1\frac{dydx}{(x+y-xy)^2}+\int_{\frac{1-z}{z}}^1\int_{1/2}^1\frac{dydx}{(x+y-xy)^2}=1-z+\log(2z).\]
In particular, $C=\bar\mu(H_1)=\bar\mu(S_1\cap H_1)=\log2$.
\end{proof}

\begin{remark}\label{lenstra=legendre}
Note that the limiting distribution of $\Theta_n^R(x),\ R=H_1,$ in Corollary \ref{cor.i} is linear on the (maximal) interval $[0,1/2]$; the supremum $1/2$ of this interval is called the \emph{Lenstra constant}.  The Lenstra constant coincides with the \emph{Legendre constant}, which is defined as the infimum of real numbers $c>0$ for which $\Theta(x,p/q)<c$ implies that $p/q$ is a {\sc rcf}-convergent of $x$ (that is, $p/q=u_n^R/s_n^R$ for some $n$); see Theorem \ref{Legendre} below.  Analogues of Lenstra and Legendre constants can be defined for other continued fraction algorithms.  In \cite{N10}, Nakada shows that for a large class of algorithms for which the Legendre constant exists, the Lenstra constant also exists and equals the Legendre constant. 

In each of the limiting distributions of approximation coefficients $\Theta_n^R(x)$ below, a Lenstra-type constant may be defined as the supremum of the domain of the  linear part.  However, we caution the reader that the existence of a Lenstra-type constant for $R$ does \emph{not} necessarily imply the existence of a Legendre-type constant for $R$; in fact, results of Erd\"os and of Brown and Yin imply that for any $c>1$, there is no proper inducible $R$ for which $1<\Theta(x,p/q)<c$ implies that $p/q=u_n^R/s_n^R$ for some $n$ (see \cite{BY1996}).
\end{remark}

The next corollary is a result of Bosma (\cite{B1990}) and concerns approximation coefficients of the $\lambda^\text{th}$ mediant convergents; see Example \ref{H&Vregions} and Figure \ref{MetricFigs}.ii.

\begin{cor}[Theorem 1.9 of \cite{B1990}]\label{cor.ii}
Set $R:=H_{\lambda+1},\ \lambda>0,$ so that for any $x=[0;a_1,a_2,\dots]\in(0,1)\backslash \mathbb{Q}$,
\[\{u_n^R/s_n^R\}_{n\ge 0}=\{(\lambda p_j+p_{j-1})/(\lambda q_j+q_{j-1})\ |\ \lambda<a_{j+1}\}_{j\ge 0}.\]
Then for $n>0$,
\[\frac{\lambda}{\lambda+1}< \Theta_n^R(x)<\lambda+1,\]
with the upper and lower bounds optimal, and for almost every such $x$,

\[\lim_{n\to\infty}\frac1n\#\{0\le k< n\ |\ \Theta_k^R(x)\le z\}
=\begin{cases}
\frac1C\left(\frac{\lambda+1}{\lambda}z-1+\log\frac{\lambda}{(\lambda+1)z}\right), & \frac{\lambda}{\lambda+1}\le z\le \frac{\lambda+1}{\lambda+2},\\ 
\frac1C\left(\frac{1}{\lambda(\lambda+1)}z+\log\frac{\lambda(\lambda+2)}{(\lambda+1)^2}\right), & \frac{\lambda+1}{\lambda+2}\le z\le \lambda,\\ 
\frac1C\left(1-\frac{1}{\lambda+1}z+\log\left(\frac{\lambda+2}{(\lambda+1)^2}z\right)\right), & \lambda\le z\le \lambda+1,
\end{cases}\]
where $C:=\bar\mu(H_{\lambda+1})=\log\frac{\lambda+2}{\lambda+1}$.
\end{cor}

Next we consider $R:=V_a,\ a>0$, which gives---in addition to the {\sc rcf}-convergents $p_{j-1}/q_{j-1}$ for which $a_{j+1}=a$---final, next-to-final, and so on mediant convergents when $a=1,\ a=2,$ and so on, respectively; see Example \ref{H&Vregions} and Figure \ref{MetricFigs}.iii.

\begin{cor}\label{cor.iii}
Set $R:=V_a,\ a>0,$ so that for any $x=[0;a_1,a_2,\dots]\in(0,1)\backslash \mathbb{Q}$,
\[\{u_n^R/s_n^R\}_{n\ge 0}=\{((a_{j+1}-a) p_j+p_{j-1})/((a_{j+1}-a) q_j+q_{j-1})\ |\ a_{j+1}\ge a\}_{j\ge 0}.\]
Then for $n>0$,
\[0<\Theta_n^R(x)<a+1,\]
with the upper and lower bounds optimal, and for almost every such $x$,
\[\lim_{n\to\infty}\frac1n\#\{0\le k< n\ |\ \Theta_k^R(x)\le z\}
=\begin{cases}
\frac1C\left(\frac1{a(a+1)}z\right), & 0\le z\le a,\\
\frac1C\left(1-\frac1{a+1}z+\log\left(\frac1a z\right)\right), & a\le z\le a+1,
\end{cases}\]
where $C:=\bar\mu(V_a)=\log\frac{a+1}{a}$.
\end{cor}

When $a=1$, Corollary \ref{cor.iii} above is similar to Theorem 3.2 of \cite{B1990} which investigates approximation coefficients of final mediants of $x$, i.e., of mediants of the form 
\[\frac{(a_{j+1}-1) p_j+p_{j-1}}{(a_{j+1}-1) q_j+q_{j-1}}.\]
However, when $R=V_1$, the Farey convergents $\{u_n^R/s_n^R\}_{n\ge 0}$ consist not only of these final mediants, but also of convergents of the form $p_{j-1}/q_{j-1}$ where $a_{j+1}=1$.  By removing the region $H_1$ from $R$, we omit these convergents from the set $\{u_n^R/s_n^R\}_{n\ge 0}$ and obtain the following generalisation of Bosma's result (see Figure \ref{MetricFigs}.iv):

\begin{cor}[Contains Theorem 3.2 of \cite{B1990} as a special case, namely $a=1$]\label{cor.iv}
Set $R:=V_a\backslash H_1,\ a>0$ so that for any $x=[0;a_1,a_2,\dots]\in(0,1)\backslash \mathbb{Q}$,
\[\{u_n^R/s_n^R\}_{n\ge 0}=\{((a_{j+1}-a) p_j+p_{j-1})/((a_{j+1}-a) q_j+q_{j-1})\ |\ a_{j+1}>a\}_{j\ge 0}.\]
Then for $n>0$,
\[\frac{a}{a+1}<\Theta_n^R(x)< a+1,\]
with the upper and lower bounds optimal, and for almost every such $x$,
\[\lim_{n\to\infty}\frac1n\#\{0\le k< n\ |\ \Theta_k^R(x)\le z\}
=\begin{cases}
\frac1C\left(\frac{a+1}{a}z-1+\log\frac{a}{(a+1)z}\right), & \frac{a}{a+1}\le z\le\frac{a+1}{a+2},\\
\frac1C\left(\frac1{a(a+1)}z+\log\frac{a(a+2)}{(a+1)^2}\right), & \frac{a+1}{a+2}\le z\le a,\\
\frac1C\left(1-\frac1{a+1}z+\log\left(\frac{a+2}{(a+1)^2}z\right) \right), & a\le z\le a+1,
\end{cases}\]
where $C:=\bar\mu(V_a\backslash H_1)=\log\frac{a+2}{a+1}$.
\end{cor}

Bosma (\cite{B1990}) also considers all {\sc rcf}-convergents and mediants whose approximation coefficients are no greater than some fixed $z_0\ge 0$, and of these, considers the asymptotic relative frequency of those whose approximation coefficients are no greater than $z\le z_0$ (see Figure \ref{MetricFigs}.v.):

\begin{cor}[Theorem 2.2 of \cite{B1990}]\label{cor.v}
Let $z_0\ge 0$ and set $R:=S_{z_0}$ so that for any $x=[0;a_1,a_2,\dots]\in(0,1)\backslash \mathbb{Q}$,
\[\{u_n^R/s_n^R\}_{n\ge 0}=\{u_n/s_n\ |\ \Theta(x,u_n/s_n)\le z_0\}_{n\ge 0}.\]
For almost every such $x$,
\[\lim_{n\to\infty}\frac1n\#\{0\le k< n\ |\ \Theta_k^R(x)\le z\}
=\begin{cases}
\frac{z}{z_0}, & 0\le z\le z_0\le 1,\\
\frac{z}{1+\log z_0}, & 0\le z\le 1\le z_0,\\
\frac{1+\log z}{1+\log z_0}, & 1\le z\le z_0.
\end{cases}\]
\end{cor}

\begin{remark}
Some of the cases in the frequency function of Corollary \ref{cor.v} are vacuous, depending on the value of $z_0$.  For instance, if $z_0<1$, then only the first case applies.  There are similarly vacuous cases in Corollaries \ref{cor.vi} and \ref{cor.vii} below.  
\end{remark}

The next result generalises Theorem 3.1 of \cite{B1990}, which considers approximation coefficients of final mediants corresponding to partial quotients $a_{j+1}=a$ for some fixed $a\ge 2$.  Here we consider approximation coefficients of {\sc rcf}-convergents ($\lambda=0$) or $\lambda^\text{th}$ mediants ($\lambda>0$) corresponding to partial quotients $a_{j+1}=a$ for fixed $a>0$; see Example \ref{rectangles_eg} and Figure \ref{MetricFigs}.vi.

\begin{cor}[Contains Theorem 3.1 of \cite{B1990} as a special case, namely $a\ge 2,\ \lambda=a-1$]\label{cor.vi}
Let $a>0$ and $\lambda\in\{0,\dots,a-1\}$, and set $R:=V_{a-\lambda}\cap H_{\lambda+1}$ so that for any $x=[0;a_1,a_2,\dots]\in(0,1)\backslash \mathbb{Q}$,
\[\{u_n^R/s_n^R\}_{n\ge 0}=\{(\lambda p_j+p_{j-1})/(\lambda q_j+q_{j-1})\ |\ a_{j+1}=a\}_{j\ge 0}.\]
Then for $n>0$,
\[\frac{(a-\lambda)\lambda}{a}<\Theta_n^R(x)<\frac{(a-\lambda+1)(\lambda+1)}{a+2},\]
with the upper and lower bounds optimal, and for almost every such $x$,
\begin{align*}
&\lim_{n\to\infty}\frac1n\#\left\{0\le k< n\ |\ \Theta_k^R(x)\le z\right\}\\
=&\begin{cases}
\frac1C\left(\frac{a}{(a-\lambda)\lambda}z-1+\log\frac{(a-\lambda)\lambda}{az}\right),& \frac{(a-\lambda)\lambda}{a}\le z\le\min\left\{\frac{(a-\lambda+1)\lambda}{a+1},\frac{(a-\lambda)(\lambda+1)}{a+1}\right\},\\
\frac1C\left(\frac1{\lambda(\lambda+1)}z+\log\frac{(a+1)\lambda}{a(\lambda+1)}\right),& \frac{(a-\lambda)(\lambda+1)}{a+1}\le z\le \frac{(a-\lambda+1)\lambda}{a+1},\\
\frac1C\left(\frac1{(a-\lambda)(a-\lambda+1)}z+\log\frac{(a+1)(a-\lambda)}{a(a-\lambda+1)}\right),& \frac{(a-\lambda+1)\lambda}{a+1}\le z\le \frac{(a-\lambda)(\lambda+1)}{a+1},\\
\frac1C\left(1-\frac{a+2}{(a-\lambda+1)(\lambda+1)}z+\log\left(\frac{(a+1)^2}{a(a-\lambda+1)(\lambda+1)}z\right)\right),& \max\left\{\frac{(a-\lambda+1)\lambda}{a+1},\frac{(a-\lambda)(\lambda+1)}{a+1}\right\}\le z\le\frac{(a-\lambda+1)(\lambda+1)}{a+2},
\end{cases}
\end{align*}
where $C:=\bar\mu(V_{a-\lambda}\cap H_{\lambda+1})=\log\frac{(a+1)^2}{a(a+2)}$.
\end{cor}

As Bosma based his proofs in \cite{B1990} on the natural extension of the Gauss transformation, which is a dynamical system only dealing with {\sc rcf}-convergents and no mediant convergents, his proofs are quite involved. Using essentially the same approach but to a `larger' dynamical system (induced transformations on the natural extension of the Farey tent map) makes the proofs of the results from \cite{B1990} and their generalisations easier.

Notice that for $\lambda\neq 0$, replacing $\lambda$ by $a-\lambda$ in the previous result leaves the limiting distribution unchanged.  Hence we obtain:
\begin{cor}
Let $a>1,\ \lambda\in\{1,\dots,a-1\},\ R_1=V_{a-\lambda}\cap H_{\lambda+1}$ and $R_2=V_{\lambda}\cap H_{a-\lambda+1}$ so that for any $x=[0;a_1,a_2,\dots]\in(0,1)\backslash\mathbb{Q}$,
\[\{u_n^{R_1}/s_n^{R_1}\}_{n\ge 0}=\{(\lambda p_j+p_{j-1})/(\lambda q_j+q_{j-1})\ |\ a_{j+1}=a\}_{j\ge 0}\]
and
\[\{u_n^{R_2}/s_n^{R_2}\}_{n\ge 0}=\{((a-\lambda) p_j+p_{j-1})/((a-\lambda) q_j+q_{j-1})\ |\ a_{j+1}=a\}_{j\ge 0}.\]
Then for almost every such $x$, the limiting distributions 
\[\lim_{n\to\infty}\frac1n\#\left\{0\le k< n\ |\ \Theta_k^{R_1}(x)\le z\right\}=\lim_{n\to\infty}\frac1n\#\left\{0\le k< n\ |\ \Theta_k^{R_2}(x)\le z\right\}\]
both exist and are equal to that given in Corollary \ref{cor.vi}.
\end{cor}

Lastly, we generalise Theorem 4.iii of \cite{BY1996}, which considers for fixed $k$ the approximation coefficients of convergents and the first $k$ and final $k$ mediant convergents (see also Theorems 3.1, 3.3 and 2.20 of \cite{I1989,B1990,J91}, respectively, when $k=1$).  Here---for fixed $\Lambda\ge 0$ and $A\ge 1$---we consider convergents and the first $\Lambda$ and final $A$ mediant convergents (see Figure \ref{MetricFigs}.vii):

\begin{cor}[Contains Theorem 4.iii of \cite{BY1996} as a special case, namely $\Lambda=A$]\label{cor.vii}
Let $\Lambda\ge 0$ and $A\ge 1$, and set
\[R:=\bigcup_{\lambda=0}^{\Lambda} H_{\lambda+1}\cup\bigcup_{a=1}^{A}V_a\]
so that for any $x=[0;a_1,a_2,\dots]\in(0,1)\backslash\mathbb{Q}$,
\begin{align*}
\{u_n^R/s_n^R\}_{n\ge 0}=&\{(\lambda p_j+p_{j-1})/(\lambda q_j +q_{j-1})\ |\ 0\le \lambda\le \Lambda\quad \text{and}\quad \lambda<a_{j+1}\}_{j\ge 0}\\
&\cup \{((a_{j+1}-a) p_j+p_{j-1})/((a_{j+1}-a) q_j+q_{j-1})\ |\ 1\le a\le A\quad \text{and}\quad a\le a_{j+1}\}_{j\ge 0}.
\end{align*}
Then for $n>0$,
\[0< \Theta_n^R(x)< \max\{\Lambda+1,A+1\},\]
with the upper and lower bounds optimal, and for almost every such $x$,
\[\lim_{n\to\infty}\frac1n\#\left\{0\le k< n\ |\ \Theta_k^R(x)\le z\right\}=\bar\mu_R(S_z\cap R_1)+\bar\mu_R(S_z\cap R_2),\]
where $R_1:=\bigcup_{\lambda=0}^{\Lambda} H_{\lambda+1},\ R_2=R\backslash R_1$,
\[\bar\mu_R(S_z\cap R_1)=\begin{cases}
\frac1{C}z, & 0\le z\le \frac{\Lambda+1}{\Lambda+2},\\
\frac1{C}\left(1-\frac1{\Lambda+1}z+\log\left(\frac{\Lambda+2}{\Lambda+1}z\right)\right), & \frac{\Lambda+1}{\Lambda+2}\le z\le \Lambda+1,\\
\frac1C\log(\Lambda+2), & \Lambda+1\le z,
\end{cases}\]
and
\[\bar\mu_R(S_z\cap R_2)=\begin{cases}
0, & 0\le z\le \frac{\Lambda+1}{\Lambda+2},\\
\frac1{C}\left(\frac{\Lambda+2}{\Lambda+1}z-1+\log\frac{\Lambda+1}{(\Lambda+2)z}\right), & \frac{\Lambda+1}{\Lambda+2}\le z\le \min\left\{\frac{(\Lambda+1)(A+1)}{\Lambda+A+2},1\right\},\\
\frac1{C}\left(\frac1{\Lambda+1}z+\log\frac{\Lambda+1}{\Lambda+2}\right), & 1\le z\le \frac{(\Lambda+1)(A+1)}{\Lambda+A+2},\\
\frac1{C}\left(\frac{A}{A+1}z+\log\frac{\Lambda+A+2}{(\Lambda+2)(A+1)}\right), & \frac{(\Lambda+1)(A+1)}{\Lambda+A+2}\le z\le 1,\\
\frac1{C}\left(1-\frac1{A+1}z+\log\left(\frac{\Lambda+A+2}{(\Lambda+2)(A+1)}z\right)\right), & \max\left\{\frac{(\Lambda+1)(A+1)}{\Lambda+A+2},1\right\}\le z\le A+1,\\
\frac1{C}\left(\log\frac{\Lambda+A+2}{\Lambda+2}\right), & A+1\le z,\\
\end{cases}\]
with $C=\bar\mu(R)=\log(\Lambda+A+2).$
\end{cor}

\subsection{On the theorems of Legendre, Fatou--Grace and Koksma}\label{On the theorems of Legendre, Fatou--Grace and Koksma}
Let $p/q\in\mathbb{Q}\cap [0,1]$ with $\text{gcd}\{p,q\}=1,\ q>0,$ and $x\in(0,1)\backslash\mathbb{Q}$.  We recall here the important result of Legendre in the theory of continued fractions:

\begin{thm}[Legendre, 1798 \cite{L98}]\label{Legendre}
If $\Theta(x,p/q)<1/2$, then $p/q$ is a {\sc rcf}-convergent of $x$.  Moreover, the constant $1/2$ is optimal.
\end{thm}
That the so-called \textit{Legendre constant} $1/2$ is optimal means that for any $c>1/2$, there exist $p/q$ and $x$ as above with $\Theta(x,p/q)<c$ but such that $p/q$ is not a {\sc rcf}-convergent of $x$.  As mentioned in \S\ref{Introduction} and the introduction to \S\ref{Metrical results}, Legendre's Theorem implies that the `excellent' rational approximations $p/q$ to an irrational $x$ are all {\sc rcf}-convergents of $x$.  Interpreted in a slightly different way, Legendre's Theorem gives a sufficient condition to verify that $p/q$ is a {\sc rcf}-convergent of $x$ without computing the expansion of $x$.

A similar result, first stated by Fatou and proven by Grace (and later Koksma), gives a sufficient condition to verify that $p/q$ is either a {\sc rcf}-convergent or nearest mediant (i.e., a first or final mediant) of $x$:

\begin{thm}[Fatou--Grace, 1904--1918 \cite{F04,G18,K37}]\label{Fatou--Grace}
If $\Theta(x,p/q)<1$, then $p/q$ is either a {\sc rcf}-convergent or nearest mediant of $x$.  The constant $1$ is optimal.
\end{thm}

Later, Koksma formulated and proved a similar statement regarding {\sc rcf}-convergents and first mediants:
\begin{thm}[Koksma, 1937 \cite{K37}]\label{Koksma}
If $\Theta(x,p/q)<2/3$, then $p/q$ is either a {\sc rcf}-convergent or a first mediant of $x$.  The constant $2/3$ is optimal.
\end{thm}

In \cite{BJ94}, Barbolosi and Jager give refinements of the theorems of Legendre, Fatou--Grace and Koksma.  The connections above between the map $\F$ and Farey convergents allow us to easily reobtain these refinements, \textit{assuming Fatou--Grace}.  We remark here that Barbolosi and Jager's proofs are in a sense more elementary, as they do not make this assumption.  Nevertheless, we find it worthwhile to present a new approach to these results to further highlight the versatility of the natural extension map $\F$ for studying {\sc rcf}-convergents and mediants.

We begin by setting necessary notation for the statements of the Barbolosi--Jager refinements.  Following \cite{BJ94}, for any nonzero rational\footnote{Since for any $x\in (0,1)$ both $p_{-1}/q_{-1}=1/0$ and $p_0/q_0=0/1$ are {\sc rcf}-convergents of $x$, we only consider nonzero rationals $p/q$ in what follows.} $p/q\in\mathbb{Q}\cap (0,1]$ with $\text{gcd}\{p,q\}=1$ we set $\epsilon(p/q):=(-1)^n$, where $n$ is the depth of $p/q$ (recall (\ref{two_expns}) and the definition of depth succeeding it).  For $x\in(0,1)\backslash\mathbb{Q},$ we also set
\[
\epsilon(x,p/q):=\begin{cases}
-1, & x<p/q,\\
1, & p/q<x.
\end{cases}
\]
The \textit{signature} of $p/q$ with respect to $x$ is defined as
\[\delta(x,p/q):=\epsilon(p/q)\epsilon(x,p/q)\in\{\pm 1\}.\]
The following lemma classifies the signatures of the nonzero and finite Farey convergents of $x$:

\begin{lem}\label{signature_lem}
Suppose $x\in(0,1)\backslash\mathbb{Q}$ has {\sc rcf}-expansion $[0;a_1,a_2,\dots]$, and let 
\[u_n/s_n=(\lambda_n p_{j_n}+p_{j_n-1})/(\lambda_n q_{j_n}+q_{j_n-1})\in\mathbb{Q}\backslash\{0\}\] 
be the $(n-1)^{st}$ Farey convergent of $x$, where $0\le \lambda_n<a_{j_n+1}$.  Then $\delta(x,u_n/s_n)=-1$ if and only if $u_n/s_n$ is
\begin{enumerate}
\item[(i)] a {\sc rcf}-convergent (i.e., $\lambda_n=0$) with $a_{j_n-1}=1,\ j_n>1$, or

\item[(ii)] a first mediant (i.e., $\lambda_n=1$). 
\end{enumerate}
\end{lem}
\begin{proof}
It is well-known that the odd- and even-order {\sc rcf}-convergents $(p_{2k-1}/q_{2k-1})_{k\ge 0}$ and $(p_{2k}/q_{2k})_{k\ge 0}$ form strictly decreasing and strictly increasing sequences, respectively, converging to $x$.  Moreover, for $k\ge 0$, the mediants satisfy (see \S1.4 of \cite{K97})
\[x<\frac{p_{2k+1}}{q_{2k+1}}=\frac{a_{2k+1}p_{2k}+p_{2k-1}}{a_{2k+1}q_{2k}+q_{2k-1}}<\frac{(a_{2k+1}-1)p_{2k}+p_{2k-1}}{(a_{2k+1}-1)q_{2k}+q_{2k-1}}<\dots<\frac{p_{2k}+p_{2k-1}}{q_{2k}+q_{2k-1}}<\frac{p_{2k-1}}{q_{2k-1}}\]
and
\[\frac{p_{2k}}{q_{2k}}<\frac{p_{2k+1}+p_{2k}}{q_{2k+1}+q_{2k}}<\dots<\frac{(a_{2k+2}-1)p_{2k+1}+p_{2k}}{(a_{2k+2}-1)q_{2k+1}+q_{2k}}<\frac{a_{2k+2}p_{2k+1}+p_{2k}}{a_{2k+2}q_{2k+1}+q_{2k}}=\frac{p_{2k+2}}{q_{2k+2}}<x.\]
Thus $u_n/s_n=(\lambda_np_{j_n}+p_{j_n-1})/(\lambda_nq_{j_n}+q_{j_n-1})$ lies between $p_{j_n-1}/q_{j_n-1}$ and $p_{j_n+1}/q_{j_n+1}$, and 
\[
\epsilon(x,u_n/s_n)=\begin{cases}
-1, & \text{$j_n$ is even},\\
1, & \text{$j_n$ is odd}.
\end{cases}
\]
We now consider the value of $\epsilon(u_n/s_n)$ in cases:
\begin{enumerate}
\item[(a)] Suppose that $\lambda_n=0$ and $a_{j_n-1}=1$ (the assumption that $u_n/s_n\in\mathbb{Q}\backslash\{0\}$ implies $j_n>1$).  Then
\[u_n/s_n=p_{j_n-1}/q_{j_n-1}=[0;a_1,\dots,a_{j_n-1}]=[0;a_1,\dots,a_{j_n-2}+1],\]
so the depth of $u_n/s_n$ is $j_n-2$, and 
\[
\epsilon(u_n/s_n)=\begin{cases}
-1, & \text{$j_n$ is odd},\\
1, & \text{$j_n$ is even}.
\end{cases}
\]

\item[(b)] Suppose that $\lambda_n=0$ and $a_{j_n-1}>1$ (here again $j_n>1$).  Then
\[u_n/s_n=p_{j_n-1}/q_{j_n-1}=[0;a_1,\dots,a_{j_n-1}]=[0;a_1,\dots,a_{j_n-1}-1,1],\]
so the depth of $u_n/s_n$ is $j_n-1$, and
\[
\epsilon(u_n/s_n)=\begin{cases}
-1, & \text{$j_n$ is even},\\
1, & \text{$j_n$ is odd}.
\end{cases}
\]

\item[(c)]  Suppose that $\lambda_n=1$.  
Then
\[u_n/s_n=(p_{j_n}+p_{j_n-1})/(q_{j_n}+q_{j_n-1})=[0;a_1,\dots,a_{j_n},1]=[0;a_1,\dots,a_{j_n}+1],\]
so the depth of $u_n/s_n$ is $j_n$, and
\[
\epsilon(u_n/s_n)=\begin{cases}
-1, & \text{$j_n$ is odd},\\
1, & \text{$j_n$ is even}.
\end{cases}
\]

\item[(d)]  Lastly, if $\lambda_n>1$, then 
\[u_n/s_n=(\lambda_n p_{j_n}+p_{j_n-1})/(\lambda_n q_{j_n}+q_{j_n-1})=[0;a_1,\dots,a_{j_n},\lambda_n]=[0;a_1,\dots,a_{j_n},\lambda_n-1,1],\]
so the depth of $u_n/s_n$ is $j_n+1$, and
\[
\epsilon(u_n/s_n)=\begin{cases}
-1, & \text{$j_n$ is even},\\
1, & \text{$j_n$ is odd}.
\end{cases}
\]
\end{enumerate}
The result now follows by setting $\delta(x,u_n/s_n)=\epsilon(u_n/s_n)\epsilon(x,u_n/s_n)$.
\end{proof}

\begin{figure}[t]
\includestandalone[width=.4\textwidth]{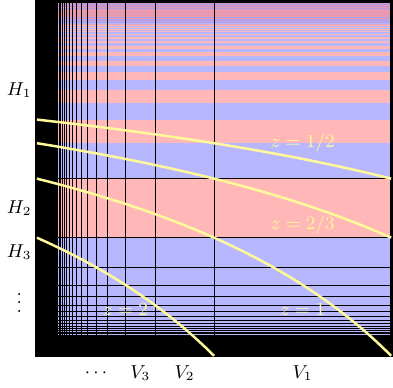}
\caption{Red and blue regions correspond to negative and positive signatures $\delta(x,u_n/s_n)=-1,1$, respectively.  The curves $h(x,y)=z$ are shown in yellow for $z\in\{1/2,2/3,1,2\}$.}
\label{BarbJagerFig}
\end{figure}

With Lemma \ref{signature_lem}, one may decompose the domain $\Omega$ of $\F$ according to the signatures of Farey convergents; confer Figure \ref{BarbJagerFig}.  Recall (Example \ref{H&Vregions}) that for fixed $\lambda$, the horizontal region $H_{\lambda+1}\subset\Omega$ corresponds to Farey convergents of $x$ of the form $u_n/s_n=(\lambda_n p_{j_n}+p_{j_n-1})/(\lambda_n q_{j_n}+q_{j_n-1})$ with $\lambda_n=\lambda$.  Thus Lemma \ref{signature_lem} implies that the region $H_2$ corresponds to Farey convergents with negative signature $\delta(x,u_n/s_n)=-1$, and the region $\cup_{\lambda>1}H_{\lambda+1}$ corresponds to Farey convergents with positive signature $\delta(x,u_n/s_n)=1$.  The region $H_1$ corresponding to {\sc rcf}-convergents $u_n/s_n=p_{j_n-1}/q_{j_n-1}$ is further decomposed depending on the value of the partial quotient $a_{j_n-1}$ in the {\sc rcf}-expansion of $x$.  Recall from (\ref{F^n_explicit}) that if $(x_n,y_n)=\F^n(x,1)\in H_1$, then
\[(x_n,y_n)=\big([0;a_{j_n+1},a_{j_n+2},\dots],[0;1,a_{j_n},a_{j_n-1},\dots,a_1]\big).\]
The cylinder of points with {\sc rcf}-expansion beginning with $[0;1,a,b,\dots]$ is the interval
\[\big([0;1,a,b+1],[0;1,a,b]\big].\]
Thus, by Lemma \ref{signature_lem}, the regions
\[[0,1]\times \big([0;1,a,2],[0;1,a,1]\big]=[0,1]\times \left(\frac{2a+1}{2a+3},\frac{a+1}{a+2}\right],\quad a\ge 1\]
correspond to Farey convergents with negative signature $\delta(x,u_n/s_n)=-1$, while the regions
\[[0,1]\times \bigcup_{b>1} \big([0;1,a,b+1],[0;1,a,b]\big]=[0,1]\times \big([0;1,a],[0;1,a,2]\big]=[0,1]\times\left(\frac{a}{a+1},\frac{2a+1}{2a+3}\right],\quad a\ge 1.\]
correspond to Farey convergents with positive signature $\delta(x,u_n/s_n)=1$.  

Recall that Theorem \ref{Fatou--Grace} implies that if $\Theta(x,p/q)<1$, then $p/q$ is \emph{some} Farey convergent\footnote{Of course, Theorem \ref{Fatou--Grace} makes the stronger claim that $p/q$ is a {\sc rcf}-convergent or nearest mediant, but only this weaker implication is needed.} $u_n/s_n$ of $x$.  With this in mind, Barbolosi and Jager's refinements now follow easily from Proposition \ref{Theta&h} and Figure \ref{BarbJagerFig}: 

\begin{cor}[Theorem 2.2 of \cite{BJ94}]\label{BJcor1}
Let $p/q\in\mathbb{Q}\cap (0,1]$ with $\text{gcd}\{p,q\}=1$ and $x\in(0,1)\backslash\mathbb{Q}$. 
\begin{itemize}
\item When $\delta(x,p/q)=-1$, $\Theta(x,p/q)<1/2$ implies that $p/q$ is a {\sc rcf}-convergent of $x$, while $\Theta(x,p/q)>2/3$ implies that $p/q$ is not a {\sc rcf}-convergent of $x$.  

\item When $\delta(x,p/q)=1$, $\Theta(x,p/q)<2/3$ implies that $p/q$ is a {\sc rcf}-convergent of $x$, while $\Theta(x,p/q)>1$ implies that $p/q$ is not a {\sc rcf}-convergent of $x$.
\end{itemize}
All constants are optimal.
\end{cor}

\begin{cor}[Theorem 4.3 of \cite{BJ94}]\label{BJcor2}
Let $p/q\in\mathbb{Q}\cap (0,1]$ with $\text{gcd}\{p,q\}=1$ and $x\in(0,1)\backslash\mathbb{Q}$.
\begin{itemize}
\item When $\delta(x,p/q)=-1$, $\Theta(x,p/q)<1$ implies $p/q$ is a {\sc rcf}-convergent or first mediant of $x$, while $\Theta(x,p/q)>2$ implies $p/q$ is neither a {\sc rcf}-convergent nor first mediant of $x$.

\item If $\delta(x,p/q)=1$, then $p/q$ is not a first mediant of $x$.  
\end{itemize}
All constants are optimal.
\end{cor}

Notice that Corollary \ref{BJcor1} implies Legendre's Theorem (Theorem \ref{Legendre}), and Corollaries \ref{BJcor1} and \ref{BJcor2} imply Koksma's Theorem (Theorem \ref{Koksma}).

\begin{cor}[Theorem 4.7 of \cite{BJ94}]
Let $p/q\in\mathbb{Q}\cap (0,1]$ with $\text{gcd}\{p,q\}=1$ and $x\in(0,1)\backslash\mathbb{Q}$.
\begin{itemize}
\item  Suppose that $\delta(x,p/q)=-1$.  If $p/q$ is a final mediant of $x$, then it is also a first mediant of $x$.  Moreover, $\Theta(x,p/q)<2/3$ implies $p/q$ is a {\sc rcf}-convergent or a final mediant of $x$, while $\Theta(x,p/q)>1$ implies $p/q$ is neither a {\sc rcf}-convergent nor a final mediant of $x$.

\item When $\delta(x,p/q)=1$, $\Theta(x,p/q)<1$ implies that $p/q$ is a {\sc rcf}-convergent or a final mediant of $x$, while $\Theta(x,p/q)>2$ implies that $p/q$ is neither a {\sc rcf}-convergent nor a nearest mediant of $x$.
\end{itemize}
All constants are optimal.
\end{cor}

The Barbolosi--Jager refinements of Legendre, Fatou--Grace and Koksma are refinements of the assumptions on the rational $p/q$ (namely, its signature) which approximates $x$.  One could instead ask for refinements of the bounds $1/2,\ 2/3,\ 1$ which occur in Theorems \ref{Legendre}, \ref{Fatou--Grace} and \ref{Koksma}.  As the next theorem shows, such refinements give information about certain partial quotients occurring in the {\sc rcf}-expansion of $x$.  

\begin{thm}
Fix some positive integer $k$.  If $p/q\in\mathbb{Q}\cap (0,1)$ with $\text{gcd}\{p,q\}=1$, $x\in (0,1)\backslash\mathbb{Q}$ and
\[\frac{k}{k+1}\le \Theta(x,p/q)<\frac{k+1}{k+2},\]
then one of the following holds:
\begin{enumerate}
\item[(i)] $p/q$ is a {\sc rcf}-convergent of the form 
\[\frac{p}{q}=\frac{p_{j-1}}{q_{j-1}}\]
for some $j$ such that $a_j=1$ and $a_{j+1}\ge k$,

\item[(ii)] $p/q$ is a first mediant of the form
\[\frac{p}{q}=\frac{p_j+p_{j-1}}{q_j+q_{j-1}}\]
for some $j$ such that $a_{j+1}\le k+1$, or

\item[(iii)] $p/q$ is a final mediant of the form 
\[\frac{p}{q}=\frac{(a_{j+1}-1)p_j+p_{j-1}}{(a_{j+1}-1)q_j+q_{j-1}}\]
for some $j$ such that $a_{j+1}\le k+1$.
\end{enumerate}
\end{thm}
\begin{proof}
Since $\Theta(x,p/q)<1$, Theorem \ref{Fatou--Grace} implies that there exists some $n$ such that 
\[\frac{p}{q}=\frac{u_n}{s_n}=\frac{\lambda_np_{j_n}+p_{j_{n-1}}}{\lambda_nq_{j_n}+q_{j_{n-1}}}\]
is a {\sc rcf}-convergent or nearest mediant of $x$.  Moreover, Proposition \ref{Theta&h} implies that $(x_n,y_n)\in S$, where 
\[S:=\left\{(x,y)\in\Omega \ \Big|\ \frac{k}{k+1}\le h(x,y)< \frac{k+1}{k+2}\right\}\subset V_1\cup H_1\cup H_2.\]
Suppose $(x_n,y_n)\in S\cap H_1$, so that $u_n/s_n=p_{j_{n-1}}/q_{j_{n-1}}$.  The curve $h(x,y)=k/(k+1)$ passes through the points $\left(0,\frac{k+1}{2k+1}\right)$ and $\left(\frac1k,\frac12\right)$ 
(see, for instance, the curves determined by $z=1/2$ and $z=2/3$ in Figure \ref{BarbJagerFig} for $k=1$ and $k=2$, respectively).  Using this, one finds $S\cap H_1\subset [0,1/k]\times (1/2,2/3]$, which is contained in the union of $\{0,1\}\times (1/2,2/3]$ and 
\[\F(V_1\cap H_1)\cap \bigcup_{a\ge k} V_a.\]
Irrationality of $x$ implies $x_n\notin \{0,1\},$ so $(x_n,y_n)\in \F(V_1\cap H_1)\cap \bigcup_{a\ge k} V_a$.  That $(x_n,y_n)\in \F(V_1\cap H_1)$ implies $a_{j_n}=1$, and that $(x_n,y_n)\in\bigcup_{a\ge k} V_a$ implies $a_{j_n+1}\ge k$, proving case (i).  The other two cases are proven similarly, considering instead when $(x_n,y_n)$ belongs to $S\cap H_2$ and $S\cap V_1$, respectively.
%
\end{proof}

In \cite{KM52}, Kuipers and Meulenbeld give sufficient conditions to guarantee that the approximation coefficients corresponding to first and final mediants are less than $1$.  This partial converse of Theorem \ref{Fatou--Grace} is also easily obtained via $\F$:
\begin{cor}[Theorem 1 of \cite{KM52}]
Let $x=[0;a_1,a_2,\dots]\in(0,1)\backslash\mathbb{Q}$.  Suppose that $a_{j+1}\ge 2$ for some $j\ge 0$.  If $a_{j+1}\le a_j+1$, then
\[\Theta\left(x,\frac{p_j+p_{j-1}}{q_j+q_{j-1}}\right)<1,\]
while if $a_{j+1}\le a_{j+2}+1$, then
\[\Theta\left(x,\frac{(a_{j+1}-1)p_j+p_{j-1}}{(a_{j+1}-1)q_j+q_{j-1}}\right)<1.\]
\end{cor}
\begin{proof}
With notation as in the statement, let $n$ be such that $u_n/s_n=(p_j+p_{j-1})/(q_j+q_{j-1})$, and suppose $a_{j+1}\le a_j+1$.  Then 
\[(x_n,y_n)=([0;a_{j+1}-1,a_{j+2},a_{j+3},\dots],[0;2,a_j,\dots,a_1])\in V_{a_{j+1}-1}\cap H_2,\]
and we have both $x_n>1/a_{j+1}$ and---since $a_j\ge a_{j+1}-1$---also $y_n\ge [0;2,a_{j+1}-1]=(a_{j+1}-1)/(2a_{j+1}-1)$.  By Proposition \ref{Theta&h},
\[\Theta\left(x,\frac{p_j+p_{j-1}}{q_j+q_{j-1}}\right)=h(x_n,y_n)<h\left(\frac1{a_{j+1}},\frac{a_{j+1}-1}{2a_{j+1}-1}\right)=1.\]
A similar argument proves the claim for final mediants. 
\end{proof}

\subsection{Consecutive approximation coefficients}\label{Consecutive approximation coefficients}
In \cite{J86} and \cite{J86_2} a two-dimensional, ergodic dynamical system $(\Delta,\mathcal{B},\nu,\mathcal{H})$ was introduced to study consecutive approximation coefficients $\theta_n(x)$ of {\sc rcf}-convergents.  The map $\mathcal{H}$ satisfies 
\begin{equation}\label{H}
\mathcal{H}(\theta_{n-1}(x),\theta_n(x))=(\theta_n(x),\theta_{n+1}(x)),\quad n>0,
\end{equation}
and is conjugate to the natural extension map $\mathcal{G}$ of the Gauss map.  This system has been used to obtain deep insights into metrical properties of approximation coefficients (\cite{J86,J86_2,JK89,K90}).  Maps which satisfy the property analogous to (\ref{H}) have also been developed for $S$-expansions, which contain a wide class of well-studied continued fractions algorithms including the {\sc rcf}, Hurwitz' singular continued fraction, Minkowski's diagonal continued fraction and Nakada's $\alpha$-continued fractions for $1/2\le \alpha\le 1$ (\cite{K1991}).  In this subsection we introduce the analogous framework for the induced systems $(R,\mathcal{B}\cap R,\bar{\mu}_R,\F_R)$.  We briefly remark that as $S$-expansions are obtained via induced maps of the Gauss natural extension, Theorem \ref{H_1&Gauss} implies that each of the aforementioned analogues of $\mathcal{H}$ may be obtained from the current setting.  

\begin{prop}\label{Theta_{n+1}}
Let $R\subset\Omega$ be an inducible subregion and $x=[0;a_1,a_2,\dots]\in(0,1)\backslash\mathbb{Q}$. Then the $(n+1)^\text{st}$ approximation coefficient $\Theta_{n+1}^R(x)$ corresponding to $R$ may be written in terms of $(x_n^R,y_n^R),\ n\ge 0,$ as
\[\Theta_{n+1}^R(x)=\det(A)\frac{(u-sx_n^R)(u+(s-u)y_n^R)}{x_n^R+y_n^R-x_n^Ry_n^R},\]
where 
\[A=\begin{pmatrix}u & t \\ s & r\end{pmatrix}:=A_{[0,r_R(x_n^R,y_n^R)]}(x_n^R),\]
with $A_{[0,r_R(x_n^R,y_n^R)]}(x_n^R)$ and $r_R(x_n^R,y_n^R)$ defined as in (\ref{A_[0,n]}) and (\ref{hitting_time}), respectively.  
\end{prop}
\begin{proof}
From (\ref{ind_map}), we have $n\ge 0$, 
\[(x_{n+1}^R,y_{n+1}^R)=\F_R(x_n^R,y_n^R)=\left(\left(A_{[0,r_R(x_n^R,y_n^R)]}(x_n^R)\right)^{-1}\cdot x_n^R,A_{[r_R(x_n^R,y_n^R),0]}(x_n^R)\cdot y_n^R\right).\]
A computation essentially identical to that of (\ref{convmats_reversed}) gives
\[A_{[r_R(x_n^R,y_n^R),0]}(x_n^R)=\begin{pmatrix} r-t & t \\ s+r-(u+t) & u+t \end{pmatrix},\]
and thus
\[(x_{n+1}^R,y_{n+1}^R)=\left(\frac{rx_n^R-t}{-sx_n^R+u},\frac{(r-t)y_n^R+t}{(s+r-(u+t))y_n^R+u+t}\right).\]
The result is obtained from a computation using this, Proposition \ref{Theta&h} and the fact that $\det(A)=\pm 1$.
\end{proof}

By partitioning $R$ into subregions on which the left columns of the matrices $A_{[0,r_R(x,y)]}(x)$ are constant, one can define an explicit function $\psi_R:R\backslash\{(0,0)\}\to\mathbb{R}^2$ satisfying
\[\psi_R(x_n^R,y_n^R)=(\Theta_n^R(x),\Theta_{n+1}^R(x)),\quad n\ge 0,\]
and thus gain insights into consecutive approximation coefficients.  This process is demonstrated for three specific inducible subregions $R\subset\Omega$ below. 

\subsubsection{Consecutive {\sc rcf}-convergents}\label{Consecutive rcf-convergents}
\begin{figure}[t]
\includestandalone[width=.35\textwidth]{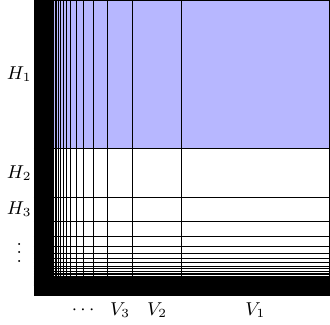}
\hspace{50pt}
\includestandalone[width=.4\textwidth]{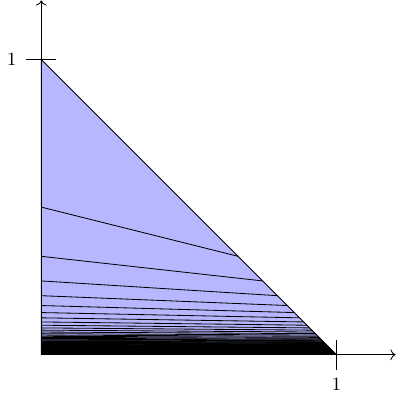}
\caption{Left: The region $R=H_1$ from \S\ref{Consecutive rcf-convergents}.  Right: The image $\psi_R(R)$.}
\label{ConsecRCFFig}
\end{figure}

Let $R=H_1$ as in Example \ref{H&Vregions}.  For any $a\ge 1$, if $(x,y)\in V_a\cap H_1$, then
\[A_{[0,r_R(x,y)]}(x)=A_0^{a-1}A_1=\begin{pmatrix}0 & 1\\ 1 & a\end{pmatrix}.\]
Define $\psi_R:R\to\mathbb{R}^2$ by
\[\psi_R(x,y)=\frac{1}{x+y-xy}(1-y,xy).\]
Propositions \ref{Theta&h} and \ref{Theta_{n+1}} then give that for any irrational $x\in(0,1)$,
\[\psi_R(x_n^R,y_n^R)=(\Theta_n^R(x),\Theta_{n+1}^R(x)),\quad n\ge 0.\]
One finds that $\psi_R$ is a diffeomorphism between the interior of $R$ and the interior of the Euclidean triangle with vertices $(0,0),\ (0,1)$ and $(1,0)$; see Figure \ref{ConsecRCFFig}. 
From this, one immediately obtains Vahlen's result that $\min\{\Theta_n^R(x),\Theta_{n+1}^R(x)\}<1/2$ for all $n\ge 0$ (\cite{V95}; see also \cite{J86_2}).  Up to the isomorphism of Theorem \ref{H_1&Gauss}, in \cite{J86_2} it is shown that $(\Delta,\mathcal{B},\nu,\mathcal{H}):=(\psi_R(R),\mathcal{B},\nu,\psi_R\circ \F_R\circ \psi_R^{-1})$ forms an ergodic system, where $d\nu=dxdy/(\log 2\sqrt{1-4xy})$.  We refer the reader to \cite{J86,J86_2,JK89,K90} for further metrical results which follow from a deeper analysis of this system.

\subsubsection{Consecutive Farey convergents}\label{Consecutive Farey convergents}
\begin{figure}[t]
\includestandalone[width=.35\textwidth]{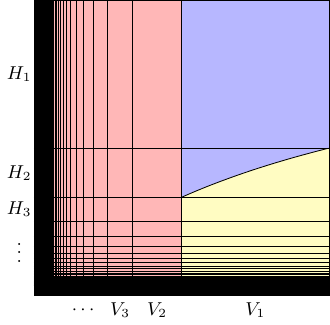}
\hspace{50pt}
\includestandalone[width=.35\textwidth]{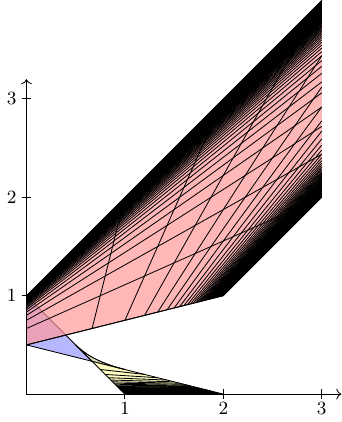}
\caption{Left: The region $R=\Omega$ from \S\ref{Consecutive Farey convergents}.  Right: The image $\psi_R(R\backslash \{(0,0)\})$.}
\label{ConsecCoeffsFig}
\end{figure}

Let $R=\Omega$.  Then for any $(x,y)\in R$, 
\[A_{[0,r_R(x,y)]}(x)=A_{\varepsilon(x)}=\begin{pmatrix}1-\varepsilon(x) & \varepsilon(x) \\ 1 & 1\end{pmatrix}.\]
Define $\psi_R:R\backslash\{(0,0)\}\to\mathbb{R}^2$ by
\begin{equation}\label{psi_R_Farey}
\psi_R(x,y)=\begin{cases}
\frac1{x+y-xy}(1-y,1-x), & x\le 1/2,\\
\frac1{x+y-xy}(1-y,xy), & x>1/2,\\
\end{cases}
\end{equation}
Now $\Theta_n^R=\Theta_n,\ x_n^R=x_n$ and $y_n^R=y_n$ for each $n\ge 0$, so Propositions \ref{Theta&h} and \ref{Theta_{n+1}} give
\begin{equation}\label{psi_R_Farey_Theta}
\psi_R(x_n,y_n)=\left(\Theta_n(x),\Theta_{n+1}(x)\right),\quad n\ge 0;
\end{equation}
see Figure \ref{ConsecCoeffsFig}.
%

The map $\psi_R$ allows one to easily compare consecutive approximation coefficients:

\begin{prop}\label{ConsecFareyConvIneqs}
Let $x=[0;a_1,a_2,\dots]\in (0,1)\backslash\mathbb{Q}$.  Then for any $j\ge 0$,
\[1<2\Theta\left(x,\frac{p_j}{q_j}\right)+\Theta\left(x,\frac{(a_{j+1}-1)p_j+p_{j-1}}{(a_{j+1}-1)q_j+q_{j-1}}\right)<2\]
and for $a_{j+1}>1$ and $0\le \lambda<a_{j+1}-1$,
\[\left|\Theta\left(x,\frac{(\lambda+1)p_j+p_{j-1}}{(\lambda+1)q_j+q_{j-1}}\right)-\Theta\left(x,\frac{(\lambda p_j+p_{j-1}}{(\lambda q_j+q_{j-1}}\right)\right|<1.\]
\end{prop}
\begin{proof}
Up to non-strict inequalities, the first statement follows from (\ref{FLseq}), (\ref{psi_R_Farey_Theta}) and the fact that the image of $V_1$ under $\psi_R$ is bounded\footnote{These lines do \emph{not} describe the boundary of $\psi_R(V_1)$ (see Figure \ref{ConsecCoeffsFig}), but rather give upper and lower bounding lines of equal slope.  More precise---though less elegant---statements can be made by analysing the boundary of the image.} by the lines $y=-x/2+1$ and $y=-x/2+1/2$, 
and the second from (\ref{FLseq}), (\ref{psi_R_Farey_Theta}) and the fact that $\psi_R(\Omega\backslash (V_1\cup\{(0,0)\}))$ is bounded by $y=x+1$ and $y=x-1$ (see Figure \ref{ConsecCoeffsFig}).  The strict inequalities follow from the irrationality of $x$.
\end{proof}

The map $\psi_R$ also gives information on the monotonicity of finite sequences of consecutive approximation coefficients.  In particular, we find that the approximation coefficients of the mediant convergents $(\lambda p_j+p_{j-1})/(\lambda q_j+q_{j-1}),\ 0< \lambda<a_{j+1}$, between $p_{j-1}/q_{j-1}$ and $p_j/q_j$ first increase and then decrease monotonically in $\lambda$, with the maximum occurring in the middle.  This statement is made precise in the following:
\begin{thm}\label{max_theta_in_middle}
Let $x=[0;a_1,a_2,\dots]\in (0,1)\backslash\mathbb{Q}$, and suppose $a_{j+1}>1$ for some fixed $j\ge 0$.  Then for each $p/q\in\{p_j/q_j\}\cup\{(\lambda p_j+p_{j-1})/(\lambda q_j+q_{j-1})\}_{0\le\lambda<a_{j+1}},\ p/q\neq 1/0,$
\[
0<\Theta(x,p/q)<\begin{cases}
\frac{a_{j+1}+2}{4}, & \text{$a_{j+1}$ even},\\
\frac{(a_{j+1}+1)(a_{j+1}+3)}{4(a_{j+1}+2)}, & \text{$a_{j+1}$ odd},
\end{cases}
\]
with the bounds optimal.  Moreover, 
\[\Theta\left(x,\frac{\lambda p_j+p_{j-1}}{\lambda q_j+q_{j-1}}\right)<\Theta\left(x,\frac{\lambda' p_j+p_{j-1}}{\lambda' q_j+q_{j-1}}\right)\quad \text{for all}\quad 0\le\lambda<\lambda'\le \lfloor a_{j+1}/2\rfloor\]
and
\[\Theta\left(x,\frac{\lambda p_j+p_{j-1}}{\lambda q_j+q_{j-1}}\right)>\Theta\left(x,\frac{\lambda' p_j+p_{j-1}}{\lambda' q_j+q_{j-1}}\right)>\Theta\left(x,\frac{p_j}{q_j}\right)\quad \text{for all}\quad \lceil a_{j+1}/2\rceil\le \lambda<\lambda'<a_{j+1}.\]
\end{thm}
\begin{proof}
Let $x$ and $a_{j+1}$ be as in the statement, $0\le \lambda\le a_{j+1}$ and let $n$ be such that $u_n/s_n=p_{j-1}/q_{j-1}$.  We begin with the latter claims.  We have (recall (\ref{FLseq}))
\[\Theta_{n+\lambda}(x)=\Theta\left(x,\frac{u_{n+\lambda}}{s_{n+\lambda}}\right)=\begin{cases}
\Theta\left(x,\frac{\lambda p_j+p_{j-1}}{\lambda q_j+q_{j-1}}\right), & \lambda<a_{j+1},\\
\Theta\left(x,\frac{p_j}{q_j}\right), & \lambda=a_{j+1},\\
\end{cases}
\]
so it suffices to show
\begin{equation}\label{Theta_eqn1}
\Theta_{n+\lambda}(x)<\Theta_{n+\lambda+1}(x),\quad 0\le\lambda< \lfloor a_{j+1}/2\rfloor
\end{equation}
and
\begin{equation}\label{Theta_eqn2}
\Theta_{n+\lambda}(x)>\Theta_{n+\lambda+1}(x),\quad \lceil a_{j+1}/2\rceil\le \lambda<a_{j+1}.
\end{equation}

For $\lambda<a_{j+1}$,
\[(x_{n+\lambda},y_{n+\lambda})\in V_{a_{j+1}-\lambda}\cap H_{\lambda+1}=\left(\frac1{a_{j+1}-\lambda+1},\frac1{a_{j+1}-\lambda}\right]\times\left(\frac1{\lambda+2},\frac1{\lambda+1}\right].\]
First, suppose $\lambda\neq a_{j+1}-1$.  Then $x_{n+\lambda}\le 1/2$, so by Equations (\ref{psi_R_Farey}) and (\ref{psi_R_Farey_Theta}), $\Theta_{n+\lambda}(x)<\Theta_{n+\lambda+1}(x)$ if $x_{n+\lambda}<y_{n+\lambda}$, and similarly with the reverse inequalities. 
If $\lambda<\lfloor a_{j+1}/2\rfloor$, then $\lambda\le a_{j+1}/2-1$, which implies that $\lambda+2\le a_{j+1}-\lambda$.  Thus in this case
\[x_{n+\lambda}\le \frac1{a_{j+1}-\lambda}\le \frac1{\lambda+2}<y_{n+\lambda},\]
proving (\ref{Theta_eqn1}).  On the other hand, if $\lceil a_{j+1}/2\rceil\le\lambda$, then $a_{j+1}/2\le \lambda$.  This implies $a_{j+1}-\lambda\le \lambda$, so in this case
\[y_{n+\lambda}\le \frac1{\lambda+1}\le\frac1{a_{j+1}-\lambda+1}<x_{n+\lambda},\]
proving (\ref{Theta_eqn2}) for $\lambda< a_{j+1}-1$.

If $\lambda=a_{j+1}-1$, then since $a_{j+1}>1$ by assumption, $(x_{n+\lambda},y_{n+\lambda})\in V_1\backslash H_1$.  By Equations (\ref{psi_R_Farey}) and (\ref{psi_R_Farey_Theta}), $\Theta_{n+\lambda}(x)>\Theta_{n+\lambda+1}(x)$ if and only if $1-y_{n+\lambda}>x_{n+\lambda}y_{n+\lambda},$ or, equivalently, $y_{n+\lambda}<1/(1+x_{n+\lambda})$, but this inequality holds since $x_{n+\lambda}< 1$ and $y_{n+\lambda}\le 1/2$.  Thus (\ref{Theta_eqn2}) is also true for $\lambda=a_{j+1}-1$.

The optimal bounds on $\Theta(x,p/q)$ follow from these monotonicity statements and the bounds of Corollary \ref{cor.vi}.
\end{proof}

As a corollary of the previous result, we find a lower bound on the maximum of the approximation coefficients of {\sc rcf}-convergents and mediants corresponding to particular partial quotients.

\begin{cor}\label{lower_bd_of_max_cor}
Let $x=[0;a_1,a_2,\dots]\in (0,1)\backslash\mathbb{Q}$, and suppose $a_{j+1}>1$ for some fixed $j\ge 0$.  Then 
\[
\max\left\{\Theta\left(x,p/q\right)\ \Big|\ \frac{p}{q}\in\left\{\frac{p_j}{q_j}\right\}\cup\left\{\frac{\lambda p_j+p_{j-1}}{\lambda q_j+q_{j-1}}\right\}_{0\le\lambda<a_{j+1}},\ p/q\neq 1/0\right\}
>\begin{cases}
\frac{a_{j+1}}{4}, & a_{j+1}\ \text{even},\\
\noalign{\vskip8pt}
\frac{a_{j+1}^2-1}{4a_{j+1}}, & a_{j+1}\ \text{odd},
\end{cases}
\]
with the lower bound optimal. 
\end{cor}
\begin{proof}
By Theorem \ref{max_theta_in_middle}, the maximum occurs when $p/q=(\lambda p_j+p_{j-1})/(\lambda q_j+q_{j-1})$ with $\lambda=\lfloor a_{j+1}/2 \rfloor$ or $\lambda=\lceil a_{j+1}/2 \rceil$.  By Corollary \ref{cor.vi}, 
\begin{equation}\label{lower_bd_of_max}
\Theta\left(x,\frac{\lambda p_j+p_{j-1}}{\lambda q_j+q_{j-1}}\right)>\frac{(a_{j+1}-\lambda)\lambda}{a_{j+1}},
\end{equation}
with the lower bound optimal.  When $a_{j+1}$ is even, $\lfloor a_{j+1}/2 \rfloor=\lceil a_{j+1}/2 \rceil=a_{j+1}/2$, and for either choice of $\lambda$ the right-hand side of (\ref{lower_bd_of_max}) equals $a_{j+1}/4$.  When $a_{j+1}$ is odd, $\lfloor a_{j+1}/2 \rfloor=(a_{j+1}-1)/2$ and $\lceil a_{j+1}/2 \rceil=(a_{j+1}+1)/2$, and for each choice of $\lambda$ the right-hand side of (\ref{lower_bd_of_max}) equals $(a_{j+1}^2-1)/4a_{j+1}$.
\end{proof}

\begin{remark}
Since the partial quotients of a.e. $x\in(0,1)\backslash\mathbb{Q}$ are unbounded (see, say, Chapter V of \cite{RS92}), the previous corollary implies that the approximation coefficients of the {\sc rcf}-convergents and mediants of a.e. $x$ are unbounded.  This is also seen intuitively in Figure \ref{ConsecCoeffsFig}; when a partial quotient $a_{j+1}$ is large, points in the $\F$-orbit of $(x,1)$ will `dip' close to the origin on the left-hand side of the figure, and, consequently, the coordinates of the corresponding images under $\psi_R$ on the right-hand side become large.  
\end{remark}

From Figure \ref{ConsecCoeffsFig} it is clear that $\psi_R$ is not injective,\footnote{The image of the red and blue regions clearly overlap, but $\psi_R$ also `folds' the blue and yellow regions in $V_1\cap H_2$ over one another.  The boundary of $\psi_R(V_1\cap H_2)$ is given by the curves $y=1-x,\ y=(2-x)/4,$ and $y=1/4x$.} and thus one cannot immediately conjugate $\F_R$ with $\psi_R$ to study the dynamics of consecutive Farey convergents as done with consecutive {\sc rcf}-convergents (\S\ref{Consecutive rcf-convergents}).  The same difficulty arises when studying consecutive convergents of the nearest integer continued fraction map (\cite{JK89}).  As in \cite{JK89}, this difficulty can be overcome by introducing a third coordinate to the image of the function $\psi_R$ which `flags' the color of the subregion that a point in the domain belongs to.  This extra coordinate makes $\psi_R$ invertible and could lead to the study of further metrical results on Farey convergents.  We leave the details to future work.

\subsubsection{Consecutive {\sc rcf}-convergents and extreme mediants}\label{Consecutive convergents and extreme mediants}
\begin{figure}[t]
\includestandalone[width=.35\textwidth]{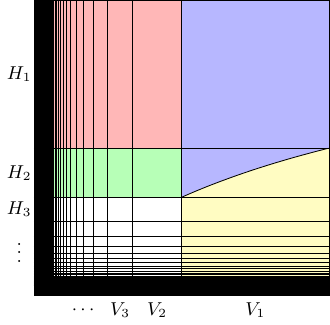}
\hspace{50pt}
\includestandalone[width=.35\textwidth]{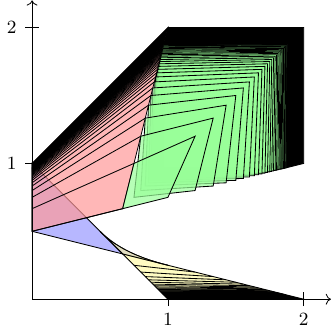}
\caption{Left: The region $R=H_1\cup H_2\cup V_1$ from \S\ref{Consecutive convergents and extreme mediants}.  Right: The image $\psi_R(R)$.}
\label{NearestMedsFig}
\end{figure}

Let $R=H_1\cup H_2\cup V_1$, which corresponds to {\sc rcf}-convergents and extreme (i.e., first and final) mediants.  For $(x,y)\in R$, we find
\[
A_{[0,r_R(x,y)]}(x)=\begin{cases}
A_1=\begin{pmatrix}0 & 1 \\ 1 & 1\end{pmatrix}, & (x,y)\in V_1,\\
A_0=\begin{pmatrix}1 & 0 \\ 1 & 1\end{pmatrix}, & (x,y)\in H_1\backslash V_1,\\
A_0^{a-1}=\begin{pmatrix}1 & 0 \\ a-1 & 1\end{pmatrix}, & (x,y)\in V_a\cap H_2,\ a>1;
\end{cases}
\]
see Figure \ref{NearestMedsFig}. Define $\psi_R:R\to \mathbb{R}^2$ by
\[
\psi_R(x,y)=\begin{cases}
\frac1{x+y-xy}\left(1-y,xy\right), & (x,y)\in V_1,\\
\frac1{x+y-xy}\left(1-y,1-x\right), & (x,y)\in H_1\backslash V_1,\\
\frac1{x+y-xy}\left(1-y,(1-(a-1)x)(1+(a-2)y)\right), & (x,y)\in V_a\cap H_2,\ a>1.
\end{cases}
\]
By Propositions \ref{Theta&h} and \ref{Theta_{n+1}}, we have for irrational $x\in (0,1)$, 
\[\psi_R(x_n^R,y_n^R)=\left(\Theta_n^R(x),\Theta_{n+1}^R(x)\right),\quad n\ge 0.\]

As with Proposition \ref{ConsecFareyConvIneqs}, the images under $\psi_R$ of subregions of $R$ yield immediate information regarding consecutive {\sc rcf}-convergents and extreme mediants.  The first statement below---which is also a corollary of Proposition \ref{ConsecFareyConvIneqs} and Theorem \ref{max_theta_in_middle}---corresponds to the image of $H_1\backslash V_1$, and the second to the image of $H_2$ (see Figure $\ref{NearestMedsFig}$; the latter image is the union of green regions and the region whose boundary is given by the curves $y=1-x,\ y=(2-x)/4,$ and $y=1/4x$).  These two images are bounded\footnote{As in the proof of Proposition \ref{ConsecFareyConvIneqs}, these lines are simply upper and lower bounding lines of equal slope and do \emph{not} necessarily describe the boundaries of the images.} by the pairs of lines $y=x+1$, $y=x$, and $y=x+1$, $y=x-1$, respectively. 

\begin{prop}
Let $x=[0;a_1,a_2,\dots]\in(0,1)\backslash\mathbb{Q}$.  For any $j\ge 0$ for which $a_{j+1}>1$,
\[0<\Theta\left(x,\frac{p_j+p_{j-1}}{q_j+q_{j-1}}\right)-\Theta\left(x,\frac{p_{j-1}}{q_{j-1}}\right)<1,\]
and
\[\left|\Theta\left(x,\frac{(a_{j+1}-1)p_j+p_{j-1}}{(a_{j+1}-1)q_j+q_{j-1}}\right)-\Theta\left(x,\frac{p_j+p_{j-1}}{q_j+q_{j-1}}\right)\right|<1.\]
\end{prop}

The map $\psi_R$ is again not injective, but a similar process to that of \cite{JK89} as mentioned in \S\ref{Consecutive Farey convergents} may be used to overcome this difficulty and thus investigate further metrical results on {\sc rcf}-convergents and extreme mediants.  We again leave the details to future work.

\subsection{A generalised L{\'e}vy-type theorem}\label{A generalised Levy-type theorem}

Recall the results (\ref{levy}) and (\ref{levy2}) due to L{\'e}vy on the growth of the denominators of {\sc rcf}-convergents and on the rate at which these convergents approach their limit $x$.  Analogues of these results have been proven for a number of continued fraction algorithms and for various subsequences of {\sc rcf}-convergents and mediants (see, e.g., \cite{I1989,J91,K1991,BY1996}).  Theorem \ref{gen_levy_thm} below generalises (\ref{levy}) and (\ref{levy2}) to subsequences of {\sc rcf}-convergents and mediants determined by proper inducible subregions $R$.  The results of \cite{I1989,J91,K1991,BY1996} then follow as special cases.

Recall from (\ref{j_lambda}) the function
\[j_n=j_n(x):=\#\{1\le k\le n\ |\ \varepsilon_k(x)=1\}\]
which counts the number of times that the $F$-orbit of an irrational $x\in(0,1)$ visits the region $(1/2,1]$ in its first $n$ steps $x,F(x),\dots,F^{n-1}(x)$.  Key to the proof of Theorem \ref{gen_levy_thm} is the simple observation that the number $j_n$ may equivalently be thought of as the number of times that the $\F$-orbit of $(x,1)$ visits $V_1\subset \Omega$ in its first $n$ steps:
\[j_n=\#\{1\le k\le n\ |\ (x_k,y_k)\in V_1\}=\sum_{k=0}^{n-1}\mathbf{1}_{V_1}(x_k,y_k).\]

\begin{thm}\label{gen_levy_thm}
For almost every $x\in(0,1)\backslash\mathbb{Q}$,
\begin{enumerate}
\item[(i)] $\lim_{n\to\infty}\frac1{n}\log s_n^R=\frac12h(\F_R)$ and 

\item[(ii)] $\lim_{n\to\infty}\frac1{n}\log \left|x-\frac{u_n^R}{s_n^R}\right|=-h(\F_R)$
\end{enumerate}
for any proper inducible $R\subset\Omega$, where $h(\F_R)=\frac{\pi^2}{6\bar\mu(R)}$ is the measure-theoretic entropy of $(R,\mathcal{B}\cap R,\bar\mu_R,\F_R)$.
\end{thm}
\begin{proof}
Let $x$ belong to the full-measure subset of $(0,1)\backslash\mathbb{Q}$ for which Theorem \ref{rel-equidist} and (\ref{levy}) both hold, and let $R\subset\Omega$ be a proper inducible subregion.  Theorem \ref{entropy_thm} gives $h(\F_R)=\frac{\pi^2}{6\bar\mu(R)}$, and by Theorem \ref{rel-equidist}, 
\begin{equation}\label{lim_jN_n}
\lim_{n\to\infty}\frac{j_{N_n}}{n}=\lim_{n\to\infty}\frac{\sum_{k=0}^{N_n-1}\mathbf{1}_{V_1}(x_k,y_k)}{\sum_{k=0}^{N_n-1}\mathbf{1}_{R}(x_k,y_k)}=\frac{\bar\mu(V_1)}{\bar\mu(R)}=\frac{\log 2}{\bar\mu(R)}, 
\end{equation}
where $N_n=N_n^R(x,1)$.  Now recall from (\ref{induced_matrices}) that $s_n^R=\lambda_{N_n}q_{j_{N_n}}+q_{j_{N_n}-1}$, so 
\begin{equation}\label{s_n^R_bounds}
q_{j_{N_n}-1}\le s_{n}^R<q_{j_{N_n}+1}.
\end{equation}
Taking logarithms and dividing by $n$, this gives
\[\frac{j_{N_n}-1}{n}\frac1{j_{N_n}-1}\log q_{j_{N_n}-1}\le \frac{1}{n}\log s_n^R<\frac{j_{N_n}+1}{n}\frac1{j_{N_n}+1}\log q_{j_{N_n}+1}.\]
Using Equations (\ref{levy}) and (\ref{lim_jN_n}), the limits as $n\to\infty$ of both the left- and right-hand sides of the previous line equal $\frac{\pi^2}{12\bar\mu(R)}$, proving (i).  

For (ii), notice that
\[\left|x-\frac{u_n^R}{s_n^R}\right|=\frac1{s_n^R(s_n^Rx_n^R+r_n^R)}\]
(see the proof of Proposition \ref{Theta&h}).  Since $x_n^R\in (0,1)$, we have 
\[\frac1{s_n^R(s_n^R+r_n^R)}<\left|x-\frac{u_n^R}{s_n^R}\right|<\frac1{s_n^Rr_n^R}.\]
Using (\ref{s_n^R_bounds}), $r_n^R=q_{j_{N_n}}$ (see (\ref{induced_matrices})) and the fact that $q_j<q_{j+1}$ for all $j$, one obtains from the previous line
\[\frac1{2q_{j_{N_n}+1}^2}<\left|x-\frac{u_n^R}{s_n^R}\right|<\frac1{q_{j_{N_n}-1}^2}.\]
Taking logarithms and dividing by $n$ gives
\[\frac{q_{j_{N_n}+1}}{n}\frac1{q_{j_{N_n}+1}}\left(-\log 2-2\log q_{j_{N_n}+1}\right)<\frac1{n}\log \left|x-\frac{u_n^R}{s_n^R}\right|<\frac{q_{j_{N_n}-1}}{n}\frac1{q_{j_{N_n}-1}}\left(-2\log q_{j_{N_n}-1}\right).\]
Again from Equations (\ref{levy}) and (\ref{lim_jN_n}), the limits as $n\to\infty$ of both the left- and right-hand sides equal $-\frac{\pi^2}{6\bar\mu(R)}$, proving (ii).  
\end{proof}

\begin{remark}
As mentioned in Remark \ref{rearr_remark} in the context of limiting distributions of approximation coefficients, the historical precedent has been to study subsequences of {\sc rcf}-convergents and mediants arranged with increasing denominators, while $(s_n^R)_{n\ge 0}$ is not necessarily increasing.  Proposition \ref{rearr_2} in the appendix (\S\ref{Appendix}) implies that a number of previously considered analogues of (\ref{levy}) and (\ref{levy2}) do indeed follow from Theorem \ref{gen_levy_thm}.  In particular, setting
\[R=\bigcup_{\lambda=0}^{\Lambda} H_{\lambda+1}\cup\bigcup_{a=1}^{A}V_a\]
for fixed $\Lambda\ge 0,\ A\ge 1$, one has $\bar\mu(R)=\log(\Lambda+A+2)$ (see Corollary \ref{cor.vii}).  When $\Lambda=A=1$, Theorem \ref{gen_levy_thm} (together with Proposition \ref{rearr_2}) gives Propositions 3.1 and 3.3(ii) of \cite{I1989} and Theorem 2.11 of \cite{J91} concerning {\sc rcf}-convergents and nearest mediants.  More generally, for $\Lambda=A\ge 1$, Theorem \ref{gen_levy_thm} gives Theorem 4.i and 4.ii of \cite{BY1996} on {\sc rcf}-convergents and the first $\Lambda$ and final $\Lambda$ mediant convergents.

Now let $Q=H_1\backslash R$ with $R\subset H_1$ a proper inducible subregion, and let $S$ be the image in $\Omega$ of $Q$ under the isomorphism of Theorem \ref{H_1&Gauss}.  Then $\bar\nu_G(S)=\bar\mu_{H_1}(Q)=\bar\mu(Q)/\log2$ and 
\[\bar\mu(R)=\bar\mu(H_1)-\bar\mu(Q)=\log2(1-\bar\nu_G(S)).\]
Thus Theorem \ref{gen_levy_thm} generalises Corollary 4.15 of \cite{K1991} on $S$-expansions.  
\end{remark}

\bigskip
\section{Appendix: Rearranging (sub-)sequences of Farey convergents}\label{Appendix}
Let $x\in(0,1)\backslash \mathbb{Q}$ with {\sc rcf}-expansion $x=[0;a_1,a_2,\dots]$.  Recall from (\ref{FLseq}) that the denominators of the sequence of Farey convergents $(u_n/s_n)_{n\ge 0}$ of $x$ are given by
\begin{alignat}{3}\label{Farey_denoms}
(s_n)_{n\ge 0}=(&q_{-1} &&,q_0+q_{-1}&&,\dots,(a_1-1)q_0+q_{-1},\\
&q_{0} &&,q_1+q_0 &&,\dots,(a_2-1)q_1+q_0,\dots,\notag\\
&q_{j-1} &&,q_j+q_{j-1} &&,\dots,(a_{j+1}-1)q_j+q_{j-1},\dots).\notag
\end{alignat}
Let $\rho$ be the bijection\footnote{For $a_1>1$ there are in fact two such bijections since $s_n=s_m$ if and only if $n=m$ or (since $q_{-1}=0$) $\{n,m\}=\{1,a_1\}$.  However, in the limits considered below, the choice between these two bijections becomes irrelevant.} of non-negative integers for which these denominators are arranged in increasing order:
\begin{alignat}{3}\label{Farey_denoms_permuted}
(s_{\rho(n)})_{n\ge 0}=(q_{-1} , &q_{0} &&,q_0+q_{-1}&&,\dots,(a_1-1)q_0+q_{-1},\\
&q_{1} &&,q_1+q_0 &&,\dots,(a_2-1)q_1+q_0,\dots,\notag\\
&q_j &&,q_j+q_{j-1} &&,\dots,(a_{j+1}-1)q_j+q_{j-1},\dots).\notag
\end{alignat}
For an inducible subregion $R\subset\Omega$, let $\rho_R$ be the bijection of non-negative integers which $(u_{\rho_R(n)}^R/s_{\rho_R(n)}^R)_{n\ge 0}$ forms a subsequence of $(u_{\rho(n)}/s_{\rho(n)})_{n\ge 0}$.  That is, $\rho_R$ permutes elements of the sequence of $(u_n^R/s_n^R)_{n\ge 0}$ of Farey convergents determined by $R$ so as to have increasing denominators.  We aim to prove the following two results:

\begin{prop}\label{rearr_1}
Let $R\subset\Omega$ be a proper inducible subregion and $z\in[0,\infty)$.  Then for any $x\in (0,1)\backslash\mathbb{Q}$,
\[\lim_{n\to\infty}\frac1n\#\{0\le k< n\ |\ \Theta_{\rho_R(k)}^R(x)\le z\}=\lim_{n\to\infty}\frac1n\#\{0\le k< n\ |\ \Theta_k^R(x)\le z\}\]
when either limit exists.
\end{prop}

\begin{prop}\label{rearr_2}
Let $R\subset\Omega$ be a proper inducible subregion.  Then for any $x\in (0,1)\backslash\mathbb{Q}$,
\begin{enumerate}
\item[(i)] $\lim\limits_{n\to\infty}\frac1{n}\log s_{\rho_R(n)}^R=\lim\limits_{n\to\infty}\frac1{n}\log s_n^R$ and

\item[(ii)] $\lim\limits_{n\to\infty}\frac1{n}\log \left|x-\frac{u_{\rho_R(n)}^R}{s_{\rho_R(n)}^R}\right|=\lim\limits_{n\to\infty}\frac1{n}\log \left|x-\frac{u_n^R}{s_n^R}\right|$
\end{enumerate}
when any of the limits exist.
\end{prop}

To prove these, we need the following:

\begin{lem}\label{rearr_lem}
For any $x=[0;a_1,a_2,\dots]\in (0,1)\backslash \mathbb{Q}$, any inducible $R\subset\Omega$ and any $n\ge 0$, the cardinality of the symmetric difference between the sets $\{N_k\}_{k=0}^n$ and $\{N_{\rho_R(k)}\}_{k=0}^n$ is at most two, and $|j_{N_{\rho_R(n)}}-j_{N_n}|\le 1$.
\end{lem}

\begin{proof}
We begin with preliminary notation and observations.  Set $A_0:=0$ and for $j\ge 1$ set $A_j:=\sum_{k=1}^j a_k$.  From (\ref{Farey_denoms}) and (\ref{Farey_denoms_permuted}), one finds that $\rho$ fixes $A_0=0$ and, for $j\ge0$,  acts as a `cyclic permutation' on subsequent blocks of length $a_{j+1}$, namely
\begin{equation}\label{local_perm}
(A_j+1,A_j+2,\dots,A_j+a_{j+1}-1,A_{j+1})\xmapsto{\rho}(A_{j+1},A_j+1,A_j+2,\dots,A_j+a_{j+1}-1),
\end{equation}
where $\rho$ is applied entry-wise.  In particular, for $N>0$ with $A_j<N\le A_{j+1}$,
\begin{equation}\label{rho_0-n}
\rho(\{0,1,\dots,N\})=\{0,1,\dots,N-1\}\cup\{A_{j+1}\}.
\end{equation}
Also observe that since $s_k^R=s_{N_k}$ and $s_{\rho_R(k)}^R=s_{N_{\rho_R(k)}}$, the numbers $N_k$ and $N_{\rho_R(k)}$ are the indices of the $(k+1)^\text{st}$ elements of the set $\{s_n^R\}_{n\ge0}=\{s_{\rho_R(n)}^R\}_{n\ge0}$ to appear in the sequences $(s_n)_{n\ge 0}$ and $(s_{\rho(n)})_{n\ge 0}$, respectively.  

If $n=0$, then since $N_0=N_{\rho_R(0)}=0$, the statement of the lemma holds.  Now fix $n>0$ and let $j$ be such that $A_j<N_n\le A_{j+1}$.  We consider two cases:
\begin{enumerate}
\item[(i)] If $s_{A_{j+1}}\notin\{s_k^R\}_{k\ge 0}$, then (\ref{rho_0-n}) implies that $\{s_{N_{\rho_R(k)}}\}_{k=0}^n=\{s_{N_k}\}_{k=0}^n$, so\footnote{Recall that $s_n=s_m$ if and only if $n=m$ or $\{n,m\}=\{1,a_1\}$.  Throughout, we consider $s_1$ and $s_{a_1}$ for $a_1>1$ as \textit{distinct} elements of $\{s_n\}_{n\ge 0}$, even though as integers these both equal $q_0=1$.} $\{N_k\}_{k=0}^n=\{N_{\rho_R(k)}\}_{k=0}^n$.  Moreover, (\ref{local_perm}) implies that $N_{\rho_R(n)}=N_n$, and thus the claim holds in this case.

\item[(ii)] Suppose $s_{A_{j+1}}\in\{s_k^R\}_{k\ge 0}$.  Then (\ref{rho_0-n}) gives that $\{s_{N_{\rho_R(k)}}\}_{k=0}^n=\{s_{N_k}\}_{k=0}^{n-1}\cup\{s_{A_{j+1}}\}$.  In particular, the symmetric difference between $\{N_k\}_{k=0}^n$ and $\{N_{\rho_R(k)}\}_{k=0}^n$ is zero if $N_n=A_{j+1}$ and two if $N_n<A_{j+1}$.  Moreover, (\ref{local_perm}) implies that
\[N_{\rho_R(n)}=\begin{cases}
N_{n-1} & \text{if $A_j<N_{n-1}$},\\
A_{j+1} & \text{otherwise}.
\end{cases}\]
If $N_{\rho_R(n)}=N_{n-1}$, then $A_j<N_{n-1}<N_n\le A_{j+1}$ implies $j_{N_{\rho_R(n)}}=j$, while if $N_{\rho_R(n)}=A_{j+1}$ we have $j_{N_{\rho_R(n)}}=j+1$.  But also $A_j<N_n\le A_{j+1}$ implies $j_{N_n}=j$ or $j_{N_n}=j+1$ (with the latter occurring if and only if $N_n=A_{j+1}$).  Thus $|j_{N_{\rho_R(n)}}-j_{N_n}|\le 1.$
\end{enumerate}
\end{proof}

Propositions \ref{rearr_1} and \ref{rearr_2} now follow almost immediately from Lemma \ref{rearr_lem}:

\begin{proof}[Proof of Proposition \ref{rearr_1}]
By definition (see (\ref{induced_matrices}) and (\ref{Theta_n^R})), $\Theta_{\rho_R(k)}^R(x)=\Theta(x,u_{N_{\rho_R(k)}}/s_{N_{\rho_R(k)}})$ and $\Theta_k^R(x)=\Theta(x,u_{N_k}/s_{N_k})$, so Lemma \ref{rearr_lem} implies that for any $n>0$,
\[\left|\#\{0\le k< n\ |\ \Theta_{\rho_R(k)}^R(x)\le z\}-\#\{0\le k< n\ |\ \Theta_k^R(x)\le z\}\right|\le 2.\]
Dividing by $n$ and taking $n\to\infty$ gives the result.  
\end{proof}

\begin{proof}[Proof of Proposition \ref{rearr_2}]
By Lemma \ref{rearr_lem},
\[\lim_{n\to\infty}\frac{j_{N_{\rho_R(n)}}}{n}=\lim_{n\to\infty}\frac{j_{N_n}}{n}\]
whenever either limit exists.  After equation (\ref{lim_jN_n}), the proof of Theorem \ref{gen_levy_thm} goes through with each of $s_n^R,\ r_n^R,\ x_n^R$ and $N_n$ replaced by $s_{\rho_R(n)}^R,\ r_{\rho_R(n)}^R,\ x_{\rho_R(n)}^R$ and $N_{\rho_R(n)}$, respectively.  
\end{proof}

\bigskip

\bibliography{bib/bib}
\bibliographystyle{acm}

\end{document}